\documentclass[numbers,webpdf,imanum]{ima-authoring-template}

\usepackage{lipsum}
\usepackage{amsfonts}
\usepackage{graphicx}
\usepackage{epstopdf}
\ifpdf
  \DeclareGraphicsExtensions{.eps,.pdf,.png,.jpg}
\else
  \DeclareGraphicsExtensions{.eps}
\fi

\usepackage{cleveref}
\usepackage{bm}
\usepackage{verbatim}
\usepackage{todonotes}
\usepackage{tikz}
\usepackage{pgfplots,pgfplotstable}
\usepgfplotslibrary{colorbrewer,groupplots}
\usepackage{caption}
\makeatletter
\renewcommand*\l@subsubsection[2]{{}{}{}} %
\makeatother
\usepackage{subcaption}
\newtheorem{theorem}{Theorem}[section]
\newtheorem{lemma}[theorem]{Lemma}
\newtheorem{corollary}[theorem]{Corollary}

\theoremstyle{definition}
\newtheorem{definition}[theorem]{Definition}

\theoremstyle{remark}
\newtheorem{remark}[theorem]{Remark}

\newtheorem{assumption}[theorem]{Assumption}

\title[Analysis and numerical analysis of the nematic Helmholtz--Korteweg equation]{Analysis and numerical analysis of the nematic Helmholtz--Korteweg equation}

\author{Patrick E.~Farrell
\address{\orgdiv{Mathematical Institute}, \orgname{University of Oxford}, \orgaddress{\street{Woodstock Rd}, \postcode{OX2 6GG}, \state{Oxford}, \country{UK}}} \address{\orgdiv{Mathematical Institute, Faculty of Mathematics and Physics}, \orgname{Charles University}, \orgaddress{\street{Sokolovská 49/83}, \postcode{186 75}, \state{Prague}, \country{Czechia}}}}

\author{Tim van Beeck
\address{\orgdiv{Institut für Numerische und Angewandte Mathematik}, \orgname{University of G\"ottingen}, \orgaddress{\street{Lotzestr. 16-18}, \postcode{37083}, \state{G\"ottingen}, \country{Germany}}}}

\author{Umberto Zerbinati
\address{\orgdiv{Mathematical Institute}, \orgname{University of Oxford}, \orgaddress{\street{Woodstock Rd}, \postcode{OX2 6GG}, \state{Oxford}, \country{UK}}}
}
\authormark{Patrick E.~Farrell, Tim van Beeck, and Umberto Zerbinati}

\corresp[*]{Corresponding author: \href{email:zerbinati@maths.ox.ac.uk}{zerbinati@maths.ox.ac.uk}}

\DOI{DOI HERE}
\copyrightyear{2021}
\vol{00}
\pubyear{2021}
\access{Advance Access Publication Date: Day Month Year}
\appnotes{Paper}
\copyrightstatement{Published by Oxford University Press on behalf of the Institute of Mathematics and its Applications. All rights reserved.}
\firstpage{1}

\usepackage{amsopn}
\usepackage{mathtools}

\newcommand{\hessian}{\mathcal{H}}
\newcommand{\abs}[1]{\lvert#1\rvert}
\newcommand{\Sch}{\mathcal{S}_h}

\let \vec \mathbf
\newcommand{\tvbchange}[1]{{{#1}}}

\newcommand{\tealchange}[1]{{{#1}}}

\newcommand{\impchange}[1]{{{#1}}}
\newcommand{\B}{\mathcal{B}}
\newcommand{\dx}{\,\mathrm{d}x}
\newcommand{\ds}{\,\mathrm{d}s}

\received{Date}{0}{Year}
\revised{Date}{0}{Year}
\accepted{Date}{0}{Year}

\begin{document}

\abstract{
	We analyse the nematic Helmholtz--Korteweg equation, a variant of the classical Helmholtz equation that describes time-harmonic wave propagation in calamitic fluids in the presence of nematic order. A prominent example is given by nematic liquid crystals, which can be modeled as nematic Korteweg fluids---that is, fluids whose stress tensor depends on density gradients and on a nematic director describing the orientation of the anisotropic molecules. These materials exhibit anisotropic acoustic properties that can be tuned by external electromagnetic fields, making them attractive for potential applications such as tunable acoustic resonators.
	We prove the existence and uniqueness of solutions to this equation in two and three dimensions for suitable (nonresonant) wave numbers and propose a convergent discretisation for its numerical solution.
	The discretisation of this problem is nontrivial as it demands high regularity and involves unfamiliar boundary conditions; we address these challenges by using high-order conforming finite elements and enforcing the boundary conditions with Nitsche's method.
	We illustrate our analysis with numerical simulations in two dimensions.
}


\maketitle

\section{Introduction}
In this work we consider the analysis and numerical analysis of the nematic Helmholtz--Korteweg equation
\begin{equation} \label{eq:nematic_helmholtz_korteweg}
\alpha \Delta^2 u + \beta \nabla \cdot \nabla \left(\vec{n}^\top (\hessian u) \vec{n}\right)- \Delta u - k^2 u = f,
\end{equation}
where $u$ represents the density perturbation of a calamitic fluid (a fluid composed of rodlike molecules) in the presence of nematic ordering (where the molecules locally align, like arrows in a quiver), $\hessian u$ is the Hessian matrix of $u$, $\vec{n}$ is a piecewise constant unit-vector field describing the local average orientation of the molecules in the fluid, and the parameters $\alpha \gg \beta \ge 0$ are material constitutive parameters known as acoustic susceptibilities \cite{VirgaNematoacoustics}.

The most prominent example of such fluids are nematic liquid crystals. These fluids exhibit long-range orientational order within particular temperature or concentration ranges. Molecules close to each other tend to align their molecular axes, giving rise to a macroscopic ordering.
Nematic liquid crystals have found widespread application in display technology and the design of novel devices due to their anisotropic acoustic, elastic, and optical properties~\cite{deGennes}. Furthermore, a remarkable feature of nematic liquid crystals is that these properties can be tuned by external electromagnetic fields.

The nematic Helmholtz--Korteweg equation, recently derived by a subset of the authors~\cite{farrell2024b}, is a novel partial differential equation modeling time-harmonic acoustic wave propagation in nematic liquid crystals.
Within this model, nematic liquid crystals are treated as Korteweg fluids with an additional anisotropic term in the Cauchy stress tensor $\sigma$ depending on the nematic director $\vec{n}$~\cite{VirgaNematoacoustics}, i.e.
\begin{equation}
  \label{eq:stress_tensor}
  \sigma = -p I - \alpha u \left(\nabla u \otimes \nabla u\right) - \beta \left(\nabla u \cdot \vec{n}\right) \nabla u \otimes \vec{n},
\end{equation}
where $p$ is the fluid pressure. The term with $\alpha$ is the Korteweg term, accounting for stresses caused by spatial variations in density \cite{Giovangigli}, while the term with $\beta$ is the anisotropic term arising due to nematic ordering.
The solutions of the nematic Helmholtz--Korteweg equation exhibit surprising phenomena not permitted by the classical Helmholtz equation (recovered by setting $\alpha = \beta = 0$), including anisotropy in the propagation of sound, penetration depth, and scattering~\cite{farrell2024b}.

%We remark that the nematic Helmholtz--Korteweg equation not only describes the acoustic propagation in nematic calamitic fluids, but also acoustic propagation in phase mixtures. In fact, when $\beta = 0$, the Cauchy stress tensor \eqref{eq:stress_tensor} reduces to the Korteweg stress tensor \cite{Korteweg} which is used to describe the elastic behavior of phase mixtures where the density acts as the phase field, separating the two phases \cite{Giovangigli}.

Let $\Omega \subset \mathbb{R}^d$, $d \in \{2, 3\}$, be a bounded convex domain with smooth boundary and let $\bm{\nu}$ be the outward-facing unit normal on $\partial \Omega$. For a given source term $f \in L^2(\Omega)$, we consider the boundary value problem
\begin{equation}
\label{eq:HKBC}
\begin{alignedat}{2}
\alpha\Delta^2 u + \beta \nabla \cdot \nabla  \left(\vec{n}^\top (\mathcal{H}u) \vec{n}\right) - \Delta u - k^2 u &= f &&\text{ in } \Omega,\\
\B u &= 0 &&\text{ on } \partial \Omega,
\end{alignedat}
\end{equation}
where $\B$ is an operator encoding the boundary conditions. In particular, $\B$ can encode sound soft, sound hard, or impedance boundary conditions.

Sound soft boundary conditions impose that the acoustic pressure vanishes on the boundary, which corresponds to
\begin{equation}\label{eq:BND:soft}
  \B u = \tvbchange{(\B_0 u, \B_1 u)} = (u, \alpha \Delta u + \beta \vec{n}^\top (\mathcal{H}u) \vec{n}).
\end{equation}
Sound hard boundary conditions impose that the normal component of the fluid velocity vanishes on the boundary, which corresponds to
\begin{equation}\label{eq:BND:hard}
  \B u = \tvbchange{(\B_0 u, \B_1 u)} = (\nabla u \cdot \bm{\nu}, \alpha \nabla (\Delta u) \cdot \bm{\nu} + \beta \nabla(\vec{n}^\top (\mathcal{H}u) \vec{n})\cdot \bm{\nu}).
\end{equation}
Lastly, impedance boundary conditions impose that the normal component of the fluid velocity is proportional to the fluid velocity, which corresponds to
\begin{equation}\label{eq:BND:impedance}
  \B u = \tvbchange{(\B_0 u, \B_1 u)} = (\nabla u \cdot \bm{\nu} - i \theta u, \alpha \nabla (\Delta u)\cdot \bm{\nu} + \beta \nabla (\vec{n}^T (\mathcal{H}u)  \vec{n}) \cdot \bm{\nu} - i \theta (\alpha \Delta u + \beta \vec{n}^\top (\mathcal{H}u)  \vec{n})),
\end{equation}
for specified impedance $\theta$.

For suitable (nonresonant) wave numbers $k$, we prove existence and uniqueness of solutions to \eqref{eq:HKBC} for \impchange{sound soft and sound hard boundary conditions as introduced above}, see \Cref{thm:Tcoercivity}. 
We then analyse a classical $H^2$-conforming finite element discretisation of these problems, yielding the first convergent numerical method for the nematic Helmholtz--Korteweg equation.

\impchange{Treating impedance boundary conditions \eqref{eq:BND:impedance} on the continuous level requires a different approach due to the higher order traces in \eqref{eq:BND:impedance}, see \Cref{rem:ImpedanceTrouble}. For this reason, we leave a detailed analysis of the impedance boundary conditions to future work. Nevertheless, the proposed discretisation also accommodates impedance boundary conditions, allowing us to illustrate their effects numerically.}

Compared to the analysis and discretisation of singularly perturbed biharmonic problems, the nematic Helmholtz--Korteweg equation presents additional challenges due to the anisotropic Korteweg term and the nonstandard boundary conditions introduced in \cite{farrell2024b}. To tackle these challenges, we use the T-coercivity technique \cite{BCS02,DCZ10,Cia12,Halla21}, where we first analyse an underlying eigenvalue problem \tvbchange{and then use the eigenbasis} to establish T-coercivity. On the discrete level, we use Nitsche's method \cite{Nit71} to enforce the boundary conditions. 

We apply the proposed numerical method to two-dimensional problems, illustrating the capabilities of the method and verifying that the nematic Helmholtz--Korteweg equation exhibits anisotropic wave propagation \cite{farrell2024b, VirgaNematoacoustics, MullenLuti}.

In particular, we show that the nematic Helmholtz--Korteweg equation correctly reproduces the experiment where anisotropic propagation speed was first observed (known as the Mullen--L\"uthi--Stephen experiment \cite{MullenLuti}), thus partially validating the idea first proposed in \cite{VirgaNematoacoustics} and further studied in \cite{farrell2024b}.
Lastly, the nematic Helmholtz--Korteweg equation is used to perform numerical simulations of acoustic resonance in a cavity filled with a nematic liquid crystal, paving the way to simulation-aided design of new acoustic devices based on nematic liquid crystals.

We remark that the nematic Helmholtz--Korteweg equation is a valid model for the acoustic propagation in nematic calamitic fluids provided that the nematic director field $\vec{n}$ is undistorted, i.e.\ $\nabla \vec{n} \equiv 0$.
This assumption is valid at small length scales; for larger length scales one must also consider the elastic effects related to distortions of the nematic director field.
%To do so a subset of the authors have recently derived a compressible inviscid variant of the Leslie--Ericksen equations \cite{farrell2024a}, starting from a kinetic theory approach captured by the general framework presented in \cite{carrillo2025}.

\section{Notation \& weak formulation}\label{sec:notationAndWF}
We denote by $X$ a Hilbert space with associated inner product $(\cdot,\cdot)_{X}$ and by $\mathcal{L}(X)$ the space of bounded linear operators acting on $X$. For a given bounded sesquilinear form $a(\cdot,\cdot)$, we denote by $A \in \mathcal{A}$ the associated linear operator.
To analyse the weak formulation of \eqref{eq:nematic_helmholtz_korteweg} introduced below, we use the concept of T-coercivity \cite{BCS02,DCZ10}, which was introduced to treat sign shifting coefficients \cite{DCZ10,BCC14,BCC18,Halla21}. It can also be used as an equivalent alternative to the inf-sup condition in a general setting, for example for the analysis of Helmholtz-like problems \cite{Cia12,vBZ24,buffaPerugia, BoffiJEE02}, Galbrun's equation \cite{HH21,HLS22,Halla23,vB23,HLvB25}, or the Stokes problem \cite{EJ25} and their respective finite element discretisations.
\begin{definition}[T-coercivity]\label{def:Tcoercivity}
  We call a sesquilinear form $a : X \times X \to \mathbb{C}$ T-coercive if there exists a bijective operator $T \in \mathcal{L}(X)$ and a constant $\gamma > 0$ such that 
  \begin{equation}
     \Re \{ a(Tu,u) \} \ge \gamma \Vert u \Vert^2_X \quad \forall u \in X.
  \end{equation}
  Equivalently, we require that $AT$ is coercive, where $A \in \mathcal{L}(X)$ is the operator associated to the sesquilinear form $a(\cdot,\cdot)$.
\end{definition}

It can be shown that T-coercivity is equivalent to the inf-sup condition \cite{Cia12} and thus, it is a necessary and sufficient condition for the well-posedness of the problem. 

To derive the weak formulation of \eqref{eq:nematic_helmholtz_korteweg} \impchange{with boundary conditions \eqref{eq:BND:soft} or \eqref{eq:BND:hard}}, \tvbchange{we define $X$ to be the space of $H^2(\Omega)$ functions that fulfill the lower order boundary condition $\B_0 u = 0$ exactly, that is}
\begin{align}\label{eq:Xspace}
  X \coloneqq \{ v \in H^2(\Omega) : \B_0 v = 0 \text{ on } \partial \Omega \}.
\end{align}
\tvbchange{Multiplying by} a test function $v \in X$ and \tvbchange{integrating by parts twice yields for the first (biharmonic term)}
\begin{align}
  \alpha &\int_{\Omega} \nabla \cdot \nabla (\Delta u) v \dx = \alpha \int_{\Omega} \Delta u \Delta v \dx - \alpha \int_{\partial \Omega} \Delta u \nabla v \cdot \bm{\nu} \ds + \alpha \int_{\partial \Omega} \nabla (\Delta u) \cdot \bm{\nu} v \ds,
\end{align}
and for the second (nematic) term, integration by parts twice yields
\begin{align}
  \nonumber \beta &\int_{\Omega} \nabla \cdot \nabla (\vec{n}^\top (\mathcal{H}u) \vec{n}) v \dx \\
  &= \beta \!  \int_{\Omega} \!\!\left(\vec{n}^\top (\mathcal{H}u) \vec{n}\right) \Delta v \dx - \beta \! \int_{\partial \Omega} \!\!\left(\vec{n}^T (\mathcal{H}u) \vec{n}\right) \nabla v \cdot \bm{\nu} \ds + \beta \! \int_{\partial \Omega} \nabla \left(\vec{n}^T (\mathcal{H}u) \vec{n}\right) \cdot \bm{\nu} v \ds, \label{eq:intByPartsNematic}
\end{align}
For the third (Laplacian) term, this procedure is standard.
The weak form of \eqref{eq:HKBC} is given by
\begin{equation}
  \label{eq:WeakForm}
  \text{find } u \in X \text{ such that } a(u,v) = (f,v)_{L^2(\Omega)} \quad \forall v \in X,
\end{equation}
\impchange{where we define the bilinear form $a(\cdot,\cdot)$ to be
\begin{equation}\label{eq:BFSoft}
  \begin{aligned}
    a(u,v) \coloneqq &\alpha (\Delta u, \Delta v)_{L^2(\Omega)}  + \beta (\vec{n}^\top (\mathcal{H} u) \vec{n},\Delta v)_{L^2(\Omega)}  + (\nabla u, \nabla v)_{L^2(\Omega)} \!-\! k^2 (u,v)_{L^2(\Omega)}.
  \end{aligned}
\end{equation}
The boundary terms stemming from integration by parts vanish due to the boundary conditions.}

%For impedance boundary conditions, we define the sesquilinear form $a(\cdot,\cdot)$ to be
%\begin{equation}\label{eq:BFImpedance}
%  \begin{aligned}
%    a(u,v) \coloneqq &\alpha (\Delta u, \Delta v)_{L^2(\Omega)} + \beta (\vec{n}^\top (\mathcal{H} u) \vec{n},\Delta v)_{L^2(\Omega)} +  (\nabla u, \nabla v)_{L^2(\Omega)} - k^2 (u,v)_{L^2(\Omega)} \\
%    &+ \tealchange{2 i \theta \, \langle T_0 u,\, \operatorname{tr}(v) \rangle_{\partial \Omega}} - i \theta (u,v)_{L^2(\partial \Omega)},
%  \end{aligned}
%\end{equation}
%\tealchange{where $T_0 u \coloneqq \alpha \Delta u + \beta \vec{n}^\top (\mathcal{H}u) \vec{n}$ denotes the Dirichlet-type biharmonic trace and $\langle \cdot,\,\cdot\rangle_{\partial \Omega}$ the $H^{-1/2}(\partial \Omega)$--$H^{1/2}(\partial \Omega)$ duality pairing, both defined precisely in \Cref{sec:appendix:traces}; the remaining boundary term $(u,v)_{L^2(\partial \Omega)}$ is the standard $L^2$-inner product, which is well-defined since $\operatorname{tr}(u),\operatorname{tr}(v) \in H^{3/2}(\partial \Omega) \subset L^2(\partial \Omega)$.}
%Several boundary terms have been substituted using the boundary conditions.
%\tvbchange{In particular, the lower order boundary condition $\mathcal{B}_0 u = 0$ is also applied to the test function $v \in X$. For example, in the impedance case, we use that $\nabla v \cdot \bm{\nu} = i \theta v$ to rewrite the second term in \eqref{eq:intByPartsNematic}.}
%An in-depth discussion arguing that the traces appearing in \eqref{eq:BFImpedance} are well-defined can be found in \Cref{sec:appendix:traces}.

\impchange{
\begin{remark}[Impedance boundary conditions]\label{rem:ImpedanceTrouble}
  Defining $X$ as in \eqref{eq:Xspace} with $\mathcal{B}_0$ coming from \eqref{eq:BND:impedance} would lead to the sesquilinear form 
  \begin{equation}\label{eq:BFImpedance}
    \begin{aligned}
      a(u,v) \coloneqq &\alpha (\Delta u, \Delta v)_{L^2(\Omega)} + \beta (\vec{n}^\top (\mathcal{H} u) \vec{n},\Delta v)_{L^2(\Omega)} +  (\nabla u, \nabla v)_{L^2(\Omega)} - k^2 (u,v)_{L^2(\Omega)} \\
      &+ \tealchange{2 i \theta \, \langle \alpha \Delta u + \beta \vec{n}^\top (\mathcal{H}u) \vec{n},\, \operatorname{tr}(v) \rangle_{\partial \Omega}} - i \theta (u,v)_{L^2(\partial \Omega)},
    \end{aligned}
  \end{equation}
  where $\langle \cdot, \cdot \rangle$ is an $H^{-1/2}(\partial \Omega)$--$H^{1/2}(\partial \Omega)$ duality pairing. Unfortunately, this term is not bounded on $H^2(\Omega) \times H^2(\Omega)$: for $u \in H^2(\Omega)$ one has $\Delta u \in L^2(\Omega)$, which does not admit a boundary trace, so the map $u \mapsto (\alpha \Delta u + \beta \vec{n}^\top (\mathcal{H}u) \vec{n})\vert_{\partial\Omega} \in H^{-1/2}(\partial\Omega)$ is not well-defined as a bounded operator on $H^2(\Omega)$. A possible alternative approach is to consider a mixed problem that allows to rewrite the boundary condition $\mathcal{B}_1 u = 0$ as a lower order boundary condition. For example, with $\beta = 0$, the problem could be rewritten with $\sigma \coloneqq - \Delta u$ as 
  \begin{align*}
    - \Delta u - \sigma = 0 \text{ in } \Omega, \qquad
    - \Delta \sigma + \sigma - k^2 u = f \text{ in } \Omega, 
  \end{align*}
  with boundary conditions $\nabla u \cdot \bm{\nu} = i \theta u$ and $\nabla \sigma \cdot \bm{\nu} = -i \theta \sigma$ on $\partial \Omega$. The analysis of this system, in particular after reintroducing the nematic term, is beyond the scope of this work.
\end{remark}
}

\section{Analysis of the nematic Helmholtz--Korteweg equation}
We denote by $\mathcal{A}$ the operator induced by the \impchange{bilinear form} defined in \eqref{eq:BFSoft}. 
To show the well-posedness of \eqref{eq:WeakForm}, we show that the operator $\mathcal{A}$ is T-coercive in the sense of \Cref{def:Tcoercivity} provided that $\beta$ is small enough.
%\todel{in the sound soft case, and weakly T-coercive and injective in the sound hard and impedance cases.
%We split the operator into a T-coercive part $A \in \mathcal{L}(X)$ and a compact part $K \coloneqq \mathcal{A} - A$, where $K \in \mathcal{L}(X)$ varies depending on the boundary conditions.
%In particular, we define $A \in \mathcal{L}(X)$ as
%\begin{equation}\label{eq:Aoperator}
%  (Au,v)_{H^2(\Omega)} \coloneqq \alpha (\Delta u, \Delta v)_{L^2(\Omega)}  +  \beta (\vec{n}^\top (\mathcal{H} u) \vec{n},\Delta v)_{L^2(\Omega)} + (\nabla u, \nabla v)_{L^2(\Omega)} - k (u,v)_{L^2(\Omega)}.
%\end{equation}
%Thus, $K \coloneqq \mathcal{A} - A$ only contains boundary terms, and $K = 0$ in the case of sound soft boundary conditions.}
%
\tvbchange{To construct a suitable T-operator, we introduce an auxiliary self-adjoint eigenvalue problem on $X$ by dropping the nematic term. Using the resulting eigenbasis, we construct a T-operator that swaps the sign of `problematic' eigenfunctions, allowing us to conclude the well-posedness of the source problem associated with the eigenvalue problem. Finally, we show that the operator $\mathcal{A}$ inherits the T-coercivity for $\beta$ sufficiently small.}

\tvbchange{
To be precise, we consider the eigenvalue problem: find $\widetilde u \in X$, $\widetilde\lambda \in \mathbb{R}$ such that 
\begin{equation}\label{eq:EVP0}
  (E_0 \widetilde u, v)_{H^2(\Omega)} \coloneqq \alpha (\Delta \widetilde u, \Delta v)_{L^2(\Omega)} + (\nabla \widetilde u, \nabla v)_{L^2(\Omega)} = \widetilde\lambda (\widetilde u, v)_{L^2(\Omega)} \quad \forall v \in X.
\end{equation}
}

\tvbchange{
\begin{lemma}
   The operator $E_0 \in \mathcal{L}(X)$ is Fredholm with index zero. 
\end{lemma}

\begin{proof}
  In the sound-soft case where $X = H^2(\Omega) \cap H^1_0(\Omega)$, we can apply the Miranda--Talenti inequality $|v|^2_{H^2(\Omega)} \le \|\Delta v\|^2_{L^2(\Omega)}$ \cite{Grisvard,MPS00} since the domain $\Omega$ is assumed to be convex. Hence, the operator $E_0 \in \mathcal{L}(X)$ is coercive, i.e. 
  \begin{align}\label{eq:EVPCoercSoft}
    (E_0 \widetilde{u}, \widetilde{u})_{H^2(\Omega)} = \alpha \Vert \Delta \widetilde{u} \Vert^2_{L^2(\Omega)} + \Vert \nabla \widetilde{u} \Vert^2_{L^2(\Omega)} \ge C \Vert \widetilde{u} \Vert^2_{H^2(\Omega)} \quad \forall \widetilde{u} \in X, 
  \end{align}
  where we further used the Poincar\'e inequality. In the sound hard case, we note that elliptic regularity  \cite[Sec.~2.3.3]{Grisvard} grants that 
  \begin{align*}
    \vert v \vert^2_{H^2(\Omega)} \le C \left( \Vert \Delta v \Vert^2_{L^2(\Omega)} + \Vert v \Vert_{H^1(\Omega)}^2 \right), \qquad \forall v \in X. 
  \end{align*}
  Thus, the operator $E_0 \in \mathcal{L}(X)$ satisfies a Gårding's inequality 
  \begin{align*}
    (E_0 \widetilde{u},\widetilde{u})_{H^2(\Omega)} \ge C \Vert \widetilde{u} \Vert_{H^2(\Omega)}^2 - C \Vert \widetilde{u} \Vert_{L^2(\Omega)}^2 \quad \forall \widetilde{u} \in X,
  \end{align*}
  where we use the compact Sobolev embedding $H^2(\Omega) \hookrightarrow L^2(\Omega)$. Thus, the claim follows. 
\end{proof}

In the sound soft case, the kernel of $E_0$ is trivial due to \eqref{eq:EVPCoercSoft} but for the sound hard case, the kernel of $E_0$ consists of constant functions. Thus, $\widetilde\lambda^{(1)} = 0$ with $\widetilde e^{(1)}$ being constant is a valid eigenpair in the sound hard case. Nevertheless, $E_0$ is self-adjoint and the solution operator associated to \eqref{eq:EVP0} is compact for all cases due to the compactness of the embedding $H^2(\Omega) \hookrightarrow L^2(\Omega)$. Thus, \eqref{eq:EVP0} admits a sequence of real eigenpairs $\{(\widetilde\lambda^{(i)}, \widetilde e^{(i)})\}_{i \in \mathbb{N}}$ of finite multiplicity \cite{Boffi2010}, and the eigenfunctions $\widetilde e^{(i)} \in X$ form an $L^2$-orthonormal basis of $X$ that is orthogonal with respect to the inner product $(E_0 \cdot,\cdot)_{H^2(\Omega)}$ as well. \tealchange{Following \cite{Cia12}, we endow $X$ with the equivalent inner product $(u,v)_X \coloneqq (E_0 u,v)_{H^2(\Omega)} + (u,v)_{L^2(\Omega)}$, which is norm-equivalent to $(\cdot,\cdot)_{H^2(\Omega)}$ and renormalise the eigenfunctions so that $\Vert \widetilde e^{(i)} \Vert_X = 1$ for all $i$, equivalently rescaling by a factor $(1+\widetilde\lambda^{(i)})^{-1/2}$.}

}

From this point on we will assume that $k^2 \not \in \{ \tealchange{\widetilde\lambda^{(i)}} \}_{i \in \mathbb{N}}$. Let $i_{\ast}$ be the largest index such that $\tealchange{\widetilde\lambda^{(i)}} < k^2$, i.e.\ $i_{\ast} \coloneqq \max \{ i \in \mathbb{N} : \tealchange{\widetilde\lambda^{(i)}} < k^2 \}$. Then, we define\footnote{For readability, we neglect the trivial case where $k^2 < \tealchange{\widetilde\lambda^{(1)}}$.}
 \begin{equation}\label{eq:Toperator}
    W \coloneqq \text{span}_{0 \le i \le i_{\ast}} \{ \tealchange{\widetilde e^{(i)}} \}, \qquad T \coloneqq \operatorname{id}_{X} - 2 P_W
 \end{equation}
 where $P_W : X \to W$ is the $L^2$-orthogonal projection onto $W$. \tvbchange{In the sound hard case, the constant eigenfunction associated with the zero eigenvalue always lies in $W$.} Since $T^2 = \operatorname{id}$, the  operator $T$ is bijective.

\begin{theorem}\label{thm:Tcoercivity}
  For $\beta$ sufficiently small, the operator $\tvbchange{\mathcal{A}} \in \mathcal{L}(X)$ associated with \eqref{eq:BFSoft} is T-coercive with respect to the $T$-operator defined in \eqref{eq:Toperator}.
\end{theorem}

\begin{proof}
  \tvbchange{
  Let $u \in X$. By construction, we can expand $u = \sum_j u^{(j)} \widetilde e^{(j)}$ as a linear combination of the eigenfunctions $\widetilde e^{(j)}$. The expression converges in $X$, and we have that
  \begin{align*}
    \Vert u \Vert^2_X = \sum_j (u^{(j)})^2, \qquad \Vert u \Vert^2_{L^2(\Omega)} = \sum_j C(\widetilde\lambda^{(j)}) (u^{(j)})^2, \qquad (E_0 u,u)_{H^2} = \sum_j C(\widetilde\lambda^{(j)}) \widetilde\lambda^{(j)} (u^{(j)})^2,
  \end{align*}
  where $C(\widetilde\lambda^{(j)}) \coloneqq 1/(1+\widetilde\lambda^{(j)})$.
  Thus, by construction of $T$, we have that
  \begin{align*}
     (E_0 Tu,u)_{H^2} - k^2(Tu,u)_{L^2(\Omega)} &= \sum_{ 0 \le i \le i_{\ast}} C(\widetilde\lambda^{(i)}) (k^2 - \widetilde\lambda^{(i)}) (u^{(i)})^2 + \sum_{i > i_{\ast}} C(\widetilde\lambda^{(i)})  (\widetilde\lambda^{(i)} - k^2) (u^{(i)})^2\\
    &\ge \gamma_0 \Vert u \Vert^2_{X},
  \end{align*}
  with $\gamma_0 \coloneqq \min_{i \in \mathbb{N}} \{ C(\widetilde\lambda^{(i)})\vert \widetilde\lambda^{(i)} - k^2 \vert \} > 0$. For the remaining nematic term, we estimate with the Cauchy--Schwarz inequality and $\Vert \vec{n} \Vert_{L^\infty} = 1$ that 
  \begin{align*}
    \left\vert \beta (\vec n^\top (\mathcal H T u) \vec n, \Delta u)_{L^2(\Omega)} \right\vert \le \beta \sqrt{d} \vert Tu \vert_{H^2(\Omega)} \vert u \vert_{H^2(\Omega)} \le \beta \sqrt{d} \Vert T \Vert_{\mathcal{L}(X)} \Vert u \Vert^2_{H^2(\Omega)}.
  \end{align*}
  In the sound soft and sound hard case, using the norm equivalence $\Vert u \Vert_X \sim \Vert u \Vert_{H^2(\Omega)}$, we can combine both estimates to obtain that
  \begin{align}\label{eq:TcoercivityEst}
    \Re \{ (\mathcal{A} T u, u )_{H^2(\Omega)} \} \ge (\gamma_0 - \beta \sqrt{d} \Vert T \Vert) \Vert u \Vert^2_{H^2(\Omega)},
  \end{align}
  and thus we obtain T-coercivity provided that $\beta < \gamma_0 / (\sqrt{d} \Vert T \Vert_{L(X)})$. This proves well-posedness in the sound soft and sound hard case. 
  }
\end{proof}

\begin{remark}[Helmholtz--Korteweg equation]
  Setting $\beta = 0$ in \eqref{eq:nematic_helmholtz_korteweg}, we obtain the Helmholtz--Korteweg equation.
  The continuous analysis presented in this section and the forthcoming discrete analysis allow for this case, yielding well-posedness of the Helmholtz--Korteweg equation and its discretisation.
  Furthermore, if we redefine $\Vert u \Vert_{X_\alpha} \coloneqq \alpha \vert u \vert_{H^2(\Omega)} + \Vert u \Vert_{H^1(\Omega)}$, we obtain continuity and stability of the sesquilinear form $a(\cdot,\cdot)$ with constants independent of $\alpha$. Thus, for $u_0 \in H^2(\Omega)$ being a weak solution of the Helmholtz equation and $u_{\alpha} \in H^2(\Omega)$ being the solution of the Helmholtz--Korteweg equation with the same data, we obtain that 
  \begin{align*}
    \Vert u_0 - u_{\alpha} \Vert_{X_\alpha} \le \gamma^{-1} \sup_{w \in X} \frac{\vert a(u_0-u_{\alpha},w)\vert}{\Vert w \Vert_{X_\alpha}} \le \alpha \gamma^{-1} \Vert \Delta u_0 \Vert_{L^2(\Omega)}.
  \end{align*}
  Thus, taking $\alpha \to 0$ this yields that $u_{\alpha} \to u_0$. With similar arguments, we can obtain a similar result for the case $\beta \to 0$.
\end{remark}

\section{Discretisation of the nematic Helmholtz--Korteweg equation}\label{sec:discrete}
In this section, we analyse the discretisation of \eqref{eq:WeakForm} using $H^2$-conforming finite elements. The analysis holds in two and three dimensions; our numerical experiments below will only be in two dimensions, since at the time of writing Firedrake \cite{Firedrake} does not yet support such elements in three dimensions. %A discussion on the use of a $C^0$ interior penalty method to circumvent this issue is presented in \cite{Tim2025}.

Conforming $H^2$ discretisations are naturally high-order, which is sometimes seen as a disadvantage; in this context it is appealing, since high-order discretisations achieve quasi-optimality for large wave numbers on much coarser meshes~\cite{Spence2025,DuWu2015,Babuska1995,Babuska1997}.
%For brevity we focus our analysis on the sound soft and the impedance boundary conditions;
%sound hard boundary conditions can be analysed \tvbchange{similarly to the latter by restricting to the real case.} 

Let $\{\mathcal{T}_h\}_{h}$ be a sequence of shape regular simplicial triangulations of $\Omega$. For polynomial degree $p \ge 3$, we define the finite element space $X_h \subset H^2(\Omega)$ as
\begin{equation*}
  X_h \coloneqq \{ v \in H^2(\Omega) : v \vert_{T} \in \mathcal{P}^p(T) \quad \forall T \in \mathcal{T}_h \},
\end{equation*}
where $\mathcal{P}^p(T)$ is the space of polynomials of degree $p$ on an element $T \in \mathcal{T}_h$. Then, we consider the discrete problem: find $u_h \in X_h$ such that 
\begin{equation}\label{eq:DiscreteProblemBeta}
  a_h(u_h,v_h) = (f,v_h)_{L^2(\Omega)} \quad \forall v_h \in X_h,
\end{equation}
\tvbchange{where the definition of the sesquilinear form is given below and depends on the boundary conditions. We denote by $\mathcal{A}_h \in \mathcal{L}(X_h)$ the corresponding bounded linear operator.
\impchange{In particular, we present discretisations for all three types of boundary conditions introduced in \eqref{eq:BND:soft}, \eqref{eq:BND:hard}, and \eqref{eq:BND:impedance}, even though the latter are not included in the continuous analysis.
}

A major challenge for the discrete analysis is that we do not want to impose the boundary condition $\mathcal{B}_0 u = 0$ directly into the discrete space as in \eqref{eq:Xspace} on the continuous level. This is motivated by the fact that imposing essential boundary conditions for $C^1$-conforming finite elements is practically very difficult, cf.~\cite{Kirby}.
To circumvent this, we use Nitsche's method \cite{Nit71}, but the analysis presented below naturally includes the case where this is achieved. 

The strategy to show that the resulting discrete problems \eqref{eq:DiscreteProblemBeta} are stable is similar to the approach presented above for the continuous level: we first consider an eigenvalue problem and then use its eigenbasis to construct an operator $T_h$ such that the sesquilinear form $a_h(\cdot,\cdot)$ is uniformly $T_h$-coercive. As at the continuous level we use a discrete analogue of the auxiliary problem~\eqref{eq:EVP0} (i.e.~with $\beta = 0$) and show that the sesquilinear form $a_h(\cdot,\cdot)$ inherits the $T_h$-coercivity for small parameters $\beta$. }

\subsection{Discrete weak formulations}
\impchange{We consider two different versions of Nitsche's method. The first allows us to impose sound soft boundary conditions, whereas the second treats sound hard and impedance boundary conditions. For the latter, we present the impedance version and note that the sound hard case follows by dropping all complex valued terms.} 
For \impchange{the sound soft case}, we define the following boundary terms:
\begin{subequations}\label{eq:NitscheTerms}
  \begin{align}
    N^{\Delta^2}_h(u_h,v_h) &\coloneqq \alpha (\nabla(\Delta u_h) \cdot \bm{\nu}, v_h)_{L^2(\partial \Omega)} + \alpha (u_h,\nabla(\Delta v_h) \cdot \bm{\nu})_{L^2(\partial \Omega)}, \label{eq:NitscheTerms:a} \\
    N_h^{\beta}(u_h,v_h) &\coloneqq \beta (\nabla(\vec{n}^\top (\mathcal{H}u_h) \vec{n}) \cdot \bm{\nu}, v_h)_{L^2(\partial \Omega)}, \label{eq:NitscheTerms:b}\\
    N^{\Delta}_h(u_h,v_h) &\coloneqq (\nabla u_h \cdot \bm{\nu}, v_h)_{L^2(\partial \Omega)} + (u_h,\nabla v_h \cdot \bm{\nu})_{L^2(\partial \Omega)}.
  \end{align}
\end{subequations}  
\tvbchange{The first term always stems from partial integration, and we further added symmetry terms to the Nitsche terms associated with the symmetric operators of the problem. Further, we define the stabilization term 
\begin{align*}
  S_h^D(u_h,v_h) \coloneqq h^{-3}(u_h,v_h)_{L^2(\partial \Omega)} + h^{-1}(u_h,v_h)_{L^2(\partial \Omega)}.
\end{align*} 

For the sound hard and impedance case, we choose to only consider the consistency terms arising from partial integration. To account for the application of the boundary condition onto test functions in the continuous setting, we define the correction term 
\begin{align*}
  \mathcal{C}_h(u_h,v_h) \coloneqq &- \alpha i \theta (\Delta u,v)_{L^2(\partial \Omega)} - \beta i \theta (\vec{n}^\top (\mathcal{H} u) \vec{n}, v)_{L^2(\partial \Omega)} \\
  &- \alpha (\Delta u, \nabla v \cdot \bm{\nu})_{L^2(\partial \Omega)} - \beta (\vec{n}^\top (\mathcal{H} u) \vec{n}, \nabla v \cdot \bm{\nu})_{L^2(\partial \Omega)}.
\end{align*}
Further, we define the stabilization term
\begin{align*}
  S_h^R(u_h,v_h) \coloneqq h^{-1} (\nabla u_h \cdot \bm{\nu} - i \theta u_h, \nabla v_h \cdot \bm{\nu} - i \theta v_h)_{L^2(\partial \Omega)}
\end{align*}
In the following, let $\eta > 0$ be a suitable stabilization parameter and $\epsilon \in \{ 0,1\}$ be a parameter that controls the choice between the sound-soft ($\epsilon = 1$) and \impchange{the sound hard or} impedance ($\epsilon = 0$) boundary conditions. We set 
\begin{align*}
  \Sch(u_h,v_h) \coloneqq \epsilon S_h^D(u_h,v_h) + (1-\epsilon) S_h^R(u_h,v_h),
\end{align*}
and altogether, we define the discrete sesquilinear form as
%
%For the sound-soft case, we define the discrete sesquilinear form as 
%\begin{equation}\label{eq:ahSoft}
%  \begin{aligned}
%  a_h(u_h,v_h) \coloneqq a(u_h,v_h) + N^{\Delta^2}_h(u_h,v_h) - N^{\Delta}_h(u_h,v_h) + N^{\beta}_h(u_h,v_h) +  \eta S_h^D(u_h,v_h),
%  \end{aligned}
%\end{equation}
%whereas for the impedance case, the sesquilinear form reads as 
%\begin{equation}\label{eq:ahImp}
%  \begin{aligned}
%  a_h(u_h,v_h) \coloneqq &\alpha (\Delta u_h, \Delta v_h)_{\Omega}  + \beta (\vec{n}^\top (\mathcal{H} u) \vec{n},\Delta v)_{L^2(\Omega)} +  (\nabla u, \nabla v)_{L^2(\Omega)} - k^2 (u,v)_{L^2(\Omega)} \\ 
%  &-\alpha (\Delta u_h, \nabla v_h \cdot \vec{n})_{\partial \Omega} - \beta (\vec{n}^\top (\mathcal{H} u_h) \vec{n}, \nabla v_h \cdot \vec{n})_{\partial \Omega} + \alpha i \theta (\Delta u_h,v_h)_{\partial \Omega} \\
%  &+ \beta i \theta (\vec{n}^\top (\mathcal{H} u_h) \vec{n}, v_h)_{\partial \Omega} - i \theta (u_h,v_h)_{\partial \Omega} +  \eta S_h^R(u_h,v_h).
%\end{aligned}
%\end{equation}
%
%Then, the discrete sesquilinear forms \eqref{eq:ahSoft} and \eqref{eq:ahImp} can be written as
\begin{align*}
  a_h(u_h,v_h) \coloneqq a(u_h,v_h) + (1-\epsilon) \mathcal{C}_h(u_h,v_h) + \epsilon \big( N^{\Delta^2}_h(u_h,v_h) - N^{\Delta}_h(u_h,v_h) + N^{\beta}_h(u_h,v_h) \big) + \eta \Sch(u_h,v_h),
\end{align*}
\impchange{where $a(\cdot,\cdot)$ is the sesquilinear form defined in \eqref{eq:BFSoft} or \eqref{eq:BFImpedance} in the impedance case. For the latter case, we consider the problematic boundary term preventing the continuous analysis, cf.~\Cref{rem:ImpedanceTrouble}, as an $(\cdot,\cdot)_{L^2(\partial \Omega)}$-inner product that is well-defined and bounded on the discrete space $X_h$. Implicitly, this enforces additional regularity on the solution.}

In preparation for the forthcoming analysis, we define the semi-norm associated with $\Sch(\cdot,\cdot)$ as
\begin{align*}
  \vert u_h \vert_{\Sch}^2 \coloneqq \Re \{ \Sch(u_h,u_h) \} = (\epsilon h^{-3} + h^{-1})\Vert u_h \Vert^2_{L^2(\partial \Omega)} + (1-\epsilon) h^{-1}\Vert \nabla u_h \cdot \bm{\nu} \Vert^2_{L^2(\partial \Omega)}
\end{align*}
for $u_h \in X_h$ and the norm
\begin{align*}
  \Vert u_h \Vert_h^2 \coloneqq \Vert \Delta u_h \Vert_{L^2(\Omega)}^2 + \Vert \nabla u_h \Vert_{L^2(\Omega)}^2 + \vert u_h \vert_{\Sch}^2, \qquad u_h \in X_h.
\end{align*}

Due to the smoothness assumptions on $\Omega$, elliptic regularity \cite[Sec.~2.3.3]{Grisvard} grants that $\vert v \vert_{H^2(\Omega)} \le C (\Vert \Delta v \Vert_{L^2(\Omega)} + \Vert v \Vert_{H^1(\Omega)})$ for all $v \in H^2(\Omega)$ with constant $C> 0$. Due to the compactness of the embedding $H^2(\Omega) \hookrightarrow H^1(\Omega)$, we may interchange the $H^1$-term with any term that guarantees injectivity on $H^2(\Omega)$ due to the Peetre--Tartar theorem \cite[Lem.~A.20]{EG_FE1}. Thus, we obtain that
\begin{align}\label{eq:DiscreteHessionControl}
  \vert v_h \vert^2_{H^2(\Omega)} \le C_{\mathcal{H}} \Vert v_h \Vert^2_{h}, \qquad v_h \in X_h.
\end{align}

\subsection{The discrete eigenvalue problem}
We consider the following discrete symmetric eigenvalue problem: find $u_h \in X_h$, $\lambda_h \in \mathbb{C}$, such that 
\begin{equation}\label{eq:DiscreteEVP}
  \begin{aligned}
    e_h(u_h,v_h) \coloneqq &\alpha (\Delta u_h, \Delta v_h)_{L^2(\Omega)} + (\nabla u_h, \nabla v_h)_{L^2(\Omega)} \\
    &+ \epsilon \big(N^{\Delta^2}_h(u_h,v_h) - N^{\Delta}_h(u_h,v_h) \big) + \eta \Sch(u_h,v_h) = \lambda_h (u_h,v_h)_{L^2(\Omega)}.
  \end{aligned}
\end{equation}

In preparation for showing uniform coercivity of $e_h(\cdot,\cdot)$, we require the following inverse trace inequalities that allow us to estimate the Nitsche terms defined in \eqref{eq:NitscheTerms}.}

\begin{lemma}\label{lem:invTraceHessian}
  There exist constants $C_1, C_2 > 0$ such that the following inverse trace inequalities hold for all $u_h \in X_h$: 
  \begin{align*}
    \Vert \nabla (\Delta u_h) \cdot \bm{\nu} \Vert^2_{L^2(\partial \Omega)} \le C_1 h^{-3} \Vert \Delta u_h \Vert^2_{L^2(\Omega)}, \qquad
    \Vert  \nabla(\vec{n}^\top (\mathcal{H}u_h) \vec{n}) \cdot \bm{\nu} \Vert^2_{L^2(\partial \Omega)} \le C_2 h^{-3} \vert u_h \vert^2_{H^2(\Omega)}.
  \end{align*}
\end{lemma}

\begin{proof}
  For the first inequality, we refer to \cite[Lem.~5.2]{BET24}. For the second, we set $\mathcal{K}_h \coloneqq \{ v_h \in X_h : \nabla (\vec{n}^\top \mathcal{H} v_h \vec{n}) = 0 \}$ such that $X_h = \mathcal{K}_h^{\perp} \oplus \mathcal{K}_h$ and consider the following eigenvalue problem on $\mathcal{K}_h^{\perp}$: find $u_h \in \mathcal{K}_h^{\perp}$, $\tilde{\lambda}_h \in \mathbb{R}$ such that for all $v_h \in \mathcal{K}_h^{\perp}$ we have
  \begin{equation}
    h (\nabla(\vec{n}^\top \mathcal{H} u_h \vec{n}) \cdot \bm{\nu},\nabla(\vec{n}^T \mathcal{H} v_h \vec{n})\cdot \bm{\nu})_{L^2(\partial \Omega)} = \tilde{\lambda}_h (\nabla(\vec{n}^T \mathcal{H} u_h \vec{n}) ,\nabla(\vec{n}^T \mathcal{H} v_h \vec{n}) )_{L^2(\Omega)}.
  \end{equation}
  The left hand-side is bounded and the right hand-side is coercive on $\mathcal{K}_h^{\perp}$, so the problem is well-posed, and the associated eigenvalues are positive and finite. Further, the min-max characterisation yields for the maximal eigenvalue $\tilde{\lambda}_{h,\max}$ that 
  \begin{equation}
    \tilde{\lambda}_{h,\max} = \sup_{v_h \in \mathcal{K}_h^{\perp} \setminus \{ 0 \}} \frac{h (\nabla(\vec{n}^\top \mathcal{H} v_h \vec{n}) \cdot \bm{\nu},\nabla(\vec{n}^T \mathcal{H} v_h \vec{n})\cdot \bm{\nu})_{L^2(\partial \Omega)}}{(\nabla(\vec{n}^T \mathcal{H} v_h \vec{n}) ,\nabla(\vec{n}^T \mathcal{H} v_h \vec{n}) )_{L^2(\Omega)} }. 
  \end{equation}
  Thus, we obtain that $\Vert \nabla( \vec{n}^\top \mathcal{H} u_h \vec{n}) \cdot \bm{\nu} \Vert_{L^2(\partial \Omega)}^2 \le C h^{-1} \Vert \nabla( \vec{n}^T \mathcal{H} u_h \vec{n}^T ) \Vert_{L^2(\Omega)}^2$. Then, we use the standard inverse inequality and that $\Vert \vec{n} \Vert^2_{L^\infty} = 1$ to estimate 
  \begin{equation*}
    h^{-1} \Vert \nabla( \vec{n}^\top \mathcal{H} u_h \vec{n} ) \Vert_{L^2(\Omega)}^2 \le C h^{-3} \Vert \mathcal{H} u_h \Vert_{L^2(\Omega)}^2 \le C h^{-3} \vert u_h \vert_{H^2(\Omega)}^2.
  \end{equation*}
\end{proof}

\tvbchange{ 
\begin{lemma}\label{lem:DiscreteEVP:Coerc}
  For $\eta > 0$ sufficiently large, the sesquilinear form $e_h(\cdot,\cdot)$ is uniformly coercive on $X_h$.
\end{lemma}

\begin{proof}
  Let $u_h \in X_h$ be arbitrary. Using the Cauchy--Schwarz inequality, the weighted Young's inequality with parameters $\xi_1, \xi_2$ and either \Cref{lem:invTraceHessian} or a standard trace inequality, we estimate for the Nitsche terms that 
  \begin{align*}
    N_h^{\Delta^2}(u_h,u_h) = 2 \alpha (\nabla ( \Delta u_h) \cdot \bm{\nu}, u_h)_{L^2(\partial \Omega)} &\le \alpha C^{\Delta^2} \left( \xi_1 \Vert \Delta u_h \Vert^2_{L^2(\Omega)} + \xi_1^{-1} \vert u_h \vert_{\Sch} \right), \\
    N_h^{\Delta}(u_h,u_h) = 2 (\nabla u_h \cdot \bm{\nu},u_h)_{L^2(\partial \Omega)} &\le C^{\Delta} \left( \xi_2 \Vert \nabla u_h \Vert^2_{L^2(\Omega)} + \xi_2^{-1} \vert u_h \vert_{\Sch} \right),
  \end{align*}
  where the respective $h$-scaling is absorbed into the semi-norm $\vert \cdot \vert_{\Sch}$.
  Thus, we have that 
  \begin{align*}
    \Re \{ e_h(u_h,u_h) \} \ge &\alpha (1 - \epsilon \xi_1 C^{\Delta^2}) \Vert \Delta u_h \Vert^2_{L^2(\Omega)} + (1- \epsilon \xi_2 C^{\Delta}) \Vert \nabla u_h \Vert^2_{L^2(\Omega)} \\
    &+ (\eta - \epsilon \alpha C^{\Delta^2} \xi^{-1} - \epsilon C^{\Delta} \xi_2^{-1}) \vert u_h \vert_{\Sch}^2 \ge \gamma \Vert u_h \Vert_{h}^2.
  \end{align*}
  Choosing the parameters $\xi_1$ and $\xi_2$ sufficiently small and the stabilization parameter $\eta > \alpha C^{\Delta^2} \xi_1^{-1} + C^{\Delta} \xi_2^{-1}$, we obtain that $e_h(u_h,u_h)$ is uniformly coercive.
\end{proof}

Let $\{ \lambda^{(i)}_h \}$ and $\{ e^{(i)}_h \}$ be the discrete eigenvalues and eigenfunctions obtained from the eigenvalue problem \eqref{eq:DiscreteEVP}. As in the continuous case, we renormalise the eigenfunctions with respect to the $\Vert \cdot \Vert_{h}$-norm. We define the discrete operator $T_h \in \mathcal{L}(X_h)$ as
\begin{equation*}
  W_h \coloneqq \text{span}_{0 \le i \le i_{\ast}} \{ e^{(i)}_h \}, \qquad i_{\ast} \coloneqq \max \{ i \in \mathbb{N} : \lambda^{(i)}_h < k^2 \}, \qquad T_h \coloneqq \operatorname{id}_{X_h} - 2 P_{W_h},
\end{equation*}
where $P_{W_h} : X_h \to W_h$ is the orthogonal projection onto $W_h$ and we are assuming that $h$ is small enough such that $\lambda^{(i_{\ast})}_h < k^2$. As before, we have that $T_h^2 = \operatorname{id}$, thus $T_h$ is bijective. Furthermore, $T_h$ is uniformly bounded with respect to the $\Vert \cdot \Vert_{h}$-norm.

\subsection{Stability}
Now, we show that for $\beta$ being sufficiently small, $\mathcal{A}_h$ is uniformly $T_h$-coercive on $X_h$. Essentially, the intuition is that $\mathcal{A}_h$ is a sufficiently small perturbation of the source problem associated with \eqref{eq:DiscreteEVP} and thus inherits its stability. In preparation, we estimate the perturbation terms as follows. 

\begin{lemma}
  For all $u_h,v_h \in X_h$, we have that 
  \begin{subequations}
    \begin{align}
      \beta (\vec{n}^{\top} (\mathcal{H} u_h) \vec{n}, \Delta v_h)_{L^2(\Omega)} &\le \beta C_{\mathcal{H}} \Vert u_h \Vert_{h} \Vert v_h \Vert_{h}, \label{eq:PertEst:a} \\
      N_h^{\beta}(u_h,v_h) &\le \beta C^{\beta} \left( \Vert u_h \Vert_h \vert v_h \vert_{\Sch} \right), \label{eq:PertEst:b} \\
      \mathcal{C}_h(u_h,v_h) &\le C_{\mathcal{C}} (\alpha + \beta) \Vert u_h \Vert_h \vert v_h \vert_{\Sch}. \label{eq:PertEst:c}
    \end{align}
  \end{subequations}
\end{lemma}

\begin{proof}
  To obtain \eqref{eq:PertEst:a}, we apply the Cauchy--Schwarz inequality, use that $\Vert \vec{n} \Vert_{L^\infty} = 1$, and use \eqref{eq:DiscreteHessionControl} to estimate the Hessian. To obtain \eqref{eq:PertEst:b} and \eqref{eq:PertEst:c}, we use in addition to \eqref{eq:DiscreteHessionControl} and the Cauchy--Schwarz inequality the result from \Cref{lem:invTraceHessian} or a standard trace inequality.
\end{proof}

\begin{theorem}\label{thm:DiscreteThCoercivity}
  For $\beta > 0$ and $h$ sufficiently small, and stabilization parameter $\eta > 0$ sufficiently large, the operator $\mathcal{A}_h \in \mathcal{L}(X_h)$ associated with \eqref{eq:DiscreteProblemBeta} is uniformly T$_h$-coercive.
\end{theorem}

\begin{proof}
  Let $\eta_e > 0$ be sufficiently large such that \Cref{lem:DiscreteEVP:Coerc} holds. Suppose that $h$ is small enough such that $\lambda^{(i_{\ast})}_h < k^2$. For all $u_h \in X_h$, we estimate using the same argument as in \Cref{thm:Tcoercivity} for $\eta > \eta_e$ that 
  \begin{align*}
    \Re \{ e_h(T_h u_h, u_h) - k^2 (T_h u_h,u_h)_{L^2(\Omega)} \} \ge C_{e} \Vert u_h \Vert_{h}^2 + (\eta - \eta_e) \vert u_h \vert_{\Sch}^2,
  \end{align*}
  where $C_e \coloneqq  \min_{i \in \mathbb{N}} C(\lambda^{(i)}_h) \{ \vert \lambda^{(i)}_h - k^2 \vert \} > 0$ is uniform in $h$, with $ C(\lambda^{(i)}_h)$ being the renormalisation constant. Since $T_h$ is bounded with respect to the $\Vert \cdot \Vert_{h}$-norm, we use \eqref{eq:PertEst:b} and a weighted Young's inequality to estimate 
  \begin{align*}
    \epsilon N_h^{\beta}(T_h u_h, u_h) \le \epsilon \beta C^{\beta} \Vert T_h \Vert_{\mathcal{L}(X_h)} \left(  \xi_1 \Vert u_h \Vert^2_h + \xi_1^{-1} \vert u_h \vert^2_{\Sch} \right)
  \end{align*}
  and similarly with \eqref{eq:PertEst:c} that 
  \begin{align*}
    (1-\epsilon) \mathcal{C}_h(u_h,v_h) \le (1-\epsilon) C_{\mathcal{C}} (\alpha + \beta) \Vert T_h \Vert_{\mathcal{L}(X_h)} \left( \xi_2 \Vert u_h \Vert_h^2 + \xi_2^{-1} \vert u_h \vert^2_{\Sch} \right).
  \end{align*}
  Together with \eqref{eq:PertEst:a}, we thus obtain that 
  \begin{align*}
    \Vert T_h \Vert_{\mathcal{L}(X_h)}^{-1} \Re \{ a_h(T_h u_h,u_h) \} \ge &\left( C_e \Vert T_h \Vert_{\mathcal{L}(X_h)} - \beta (C_{\mathcal{H}}+\epsilon C^{\beta} \xi_1 + (1-\epsilon) (\alpha + \beta) C_{\mathcal{C}} \xi_2) \right) \Vert u_h \Vert_{h}^2 \\
    &+ \left( \eta - \eta_e - \beta(\epsilon C^{\beta} \xi_1^{-1} + (1-\epsilon)(\alpha + \beta)C_{\mathcal{C}}\xi_2^{-1}) \right) \vert u_h \vert_{\Sch}^2
  \end{align*}
  Assuming that $\beta$ is small enough such that 
  \begin{align*}
    \beta < C_{e} \Vert T_h \Vert_{\mathcal{L}(X_h)} / C_{\mathcal{H}},
  \end{align*}
  we can choose $\xi_1$ and $\xi_2$ small enough such that the first constant on the right-hand side is positive. Choosing $\eta > 0$ sufficiently large, we obtain that $a_h(T_h u_h, u_h)$ is uniformly T$_h$-coercive.
\end{proof}}

\tvbchange{

\subsection{Convergence}
\impchange{For this section, we restrict to the setting of sound soft or sound hard boundary conditions.}
To show a convergence result, we assume that the continuous solution provided by \Cref{thm:Tcoercivity} has additional regularity.
\begin{assumption}\label{assumption:regularity}
  The solution $u$ of the continuous problem \eqref{eq:WeakForm} satisfies $u \in \tilde{X}$, where 
  \begin{align*}
    \tilde{X} \coloneqq \big\{ u \in X : \ &\epsilon \nabla(\Delta u) \cdot \bm{\nu} \in L^2(\partial \Omega), \epsilon \nabla (\vec{n}^\top (\mathcal{H} u) \vec{n}) \cdot \bm{n} \in L^2(\partial \Omega), \\ 
    &(1-\epsilon) \Delta u \in L^2(\partial \Omega), (1-\epsilon) \vec{n}^\top (\mathcal{H} u) \vec{n} \in L^2(\partial \Omega) \big\}.
  \end{align*}
\end{assumption}

In addition to the $\Vert \cdot \Vert_{h}$, we define the following stronger norm
\begin{align*}
  \Vert u_h \Vert_{h,\ast}^2 \coloneqq \Vert u_h \Vert^2_h &+\epsilon \left( h^3 \Vert \nabla (\Delta u_h) \cdot \bm{\nu} \Vert_{L^2(\partial \Omega)}^2 + h^3 \Vert \nabla (\vec{n}^\top (\mathcal{H} u_h) \vec{n}) \cdot \bm{\nu} \Vert_{L^2(\partial \Omega)}^2 + h \Vert \nabla u_h \cdot \bm{\nu} \Vert_{L^2(\partial \Omega)}^2 \right) \\
  &+(1-\epsilon) \left(  h \Vert \Delta u_h \Vert_{L^2(\partial \Omega)} + h \Vert \vec{n}^\top (\mathcal{H} u_h) \vec{n} \Vert_{L^2(\partial \Omega)}^2 \right).
\end{align*}
}
Applying the trace inverse estimates from \Cref{lem:invTraceHessian} and the standard trace inequality, we have that $\Vert u_h \Vert_{h,\ast} \simeq \Vert u_h \Vert_{h}$ for all $u_h \in X_h$. Furthermore, the norm is well-defined for elements $u \in \tilde{X}$, and we have that 
\begin{align*}
  a_h(w,v_h) \le C_{\text{cont}} \Vert w \Vert_{h,\ast} \Vert v_h \Vert_{h,\ast} \qquad \forall w \in \tilde{X} + X_h, v_h \in X_h. 
\end{align*}
\impchange{This continuity estimate is also valid for the discrete sesquilinear form with impedance boundary conditions.}

\begin{corollary}\label{cor:GalerkinOrthogonality}
  Let $u \in \tilde{X}$ be the solution to \eqref{eq:WeakForm} and $u_h \in X_h$ be the solution to the discrete problem \eqref{eq:DiscreteProblemBeta}. Then, it holds that
  \begin{equation*}
    a_h(u-u_h,v_h) = 0 \quad \forall v_h \in X_h.
  \end{equation*}
\end{corollary}

\tvbchange{
\begin{proof}
  Since $u \in \tilde{X}$ provides enough regularity to ensure that all boundary terms in $a_h(\cdot,\cdot)$ are well-defined, the statement follows from partial integration.
\end{proof}
}

Altogether, we have the following best approximation result. 
\begin{theorem}\label{thm:bestapprox}
  Let $u \in \tilde{X}$ be the solution to \eqref{eq:WeakForm}. Assuming that $\eta > 0$ is sufficiently large and $\beta > 0$ is sufficiently small, there exists $h_0 > 0$ such that the discrete problem \eqref{eq:DiscreteProblemBeta} has a unique solution $u_h \in X_h$ for all $h < h_0$. Further, there exists a constant $C >0$ such that
  \begin{equation*}
    \Vert u - u_h \Vert_{h} \le C \inf_{v_h \in X_h} \Vert u - v_h \Vert_{h,\ast}.
  \end{equation*}
\end{theorem}

\begin{proof}
  Due to \Cref{thm:DiscreteThCoercivity}, the operator $A_h \in L(X_h)$ is uniformly $T_h$-coercive on $X_h$ for $h < h_{\ast}$, where $h_{\ast} \coloneqq \max_{h} \{ \lambda_h^{(i_{\ast})} < k^2 \}$, and thus the discrete problem has a unique solution for all $h < h_{\ast}$.

  For $v_h \in X_h$ the triangle inequality yields that 
  \begin{equation*}
    \Vert u- u_h \Vert_{h} \le \Vert u - v_h \Vert_{h} + \Vert v_h - u_h \Vert_{h}.
  \end{equation*}
  Using the uniform $T_h$-coercivity of $a_h(\cdot,\cdot)$ and \Cref{cor:GalerkinOrthogonality}, we obtain for the second term,
  \begin{align*}
    \gamma \Vert u_h - v_h \Vert_{h}^2 &\le a_h(T_h(u_h - v_h),u_h - v_h) \le a_h(u - v_h, T_h^{\ast} (u_h - v_h)) \\
    &\le C_{\text{cont}} \Vert T_h \Vert_{\mathcal{L}(X_h)} \Vert u - v_h \Vert_{h,\ast} \Vert u_h - v_h \Vert_{h,\ast}.
  \end{align*}
  We note that $\Vert v_h \Vert_{h} \le \Vert v_h \Vert_{h,\ast}$, divide by $\Vert u_h - v_h \Vert_{h,\ast}$, and take the infimum over $v_h \in X_h$ we obtain
  \begin{equation*}
    \Vert u- u_h \Vert_{h} \le \left( 1+ \frac{C_{\text{cont}} \Vert T_h \Vert_{\mathcal{L}(X_h)}}{\gamma} \right) \inf_{v_h \in X_h} \Vert u - v_h \Vert_{h,\ast}. 
  \end{equation*}
\end{proof}

\begin{remark}[The threshold $h_0$]\label{rem:asymptoticRegime}
  The threshold $h_0$ in \Cref{thm:bestapprox} is in general not explicitly given. However, for sound soft boundary conditions, it can be characterised by the condition $\lambda^{(i_{\ast})}_h < k^2$. We could adapt a similar scheme as considered in \cite{vBZ24} to ensure that the smallness assumption on $h$ is satisfied: on a sequence of refined meshes, the criterion $\lambda^{(i_{\ast})}_h < k^2$ is checked by solving the eigenvalue problem \eqref{eq:DiscreteEVP}. If the criterion is met, the assumptions of \Cref{thm:DiscreteThCoercivity} are satisfied, otherwise the mesh is refined and the procedure is repeated. To reduce the required number of degrees of freedom, an adaptive error estimator for the first $i_{\ast}$ eigenfunctions can be used. In the case of impedance boundary conditions, this argument is not applicable. 
\end{remark}

\section{Numerical examples}
The analysis above is abstract and applies to any $C^1$ conforming finite element space. In two dimensions these include the Argyris element~\cite{Argyris} ($p \ge 5$), the Hsieh--Clough--Tocher macroelement defined on the Alfeld split of a triangular mesh~\cite{HCT} ($p \ge 3$), and the Morgan--Scott element~\cite{MorganScott,Parker} ($p = 5$).
Here we will present numerical simulations implemented using the Firedrake finite element library and ngsPETSc \cite{Firedrake,brubeck2025a, ngsPETSc}.
We test the Argyris element and the Hsieh--Clough--Tocher macroelements.

\subsection{Quasi optimal convergence}
\begin{figure}[htbp]
  \begin{subfigure}[b]{0.32\textwidth}
    \begin{tikzpicture}[scale=0.49]
      \begin{axis}[ymode=log, grid style=dashed, xlabel={refinement}, ylabel={$\lVert u-u_h\rVert_{H^2(\Omega)}$}, title={$k=10$, $\alpha=10^{-2}$}, legend pos=south west, legend cell align={left}]
        \addplot table [x=refs, y=H2, col sep=comma] {assets/convergence_10.0_0.01_0.0.csv};
        \addplot [color=blue, mark=square] table [x=refs, y=H2, col sep=comma,] {assets/convergence_HCT_3_10.0_0.01_0.0.csv};
        \addplot [mark=*, color=red] table [x=refs, y=H2, col sep=comma] {assets/convergence_10.0_0.01_0.005.csv};
        \addplot [mark=square, color=red] table [x=refs, y=H2, col sep=comma] {assets/convergence_HCT_3_10.0_0.01_0.005.csv};
        \addplot[gray,dashed,domain=1:4] {0.125*((1/16)^x)};
        \addplot[gray,dashed,domain=1:4] {8*((1/4)^x)};
        \legend{ARG $\beta=0$, HCT $\beta=0$, ARG $\beta=5\cdot 10^{-3}$, HCT $\beta=5\cdot 10^{-3}$}
      \end{axis}
    \end{tikzpicture}
  \end{subfigure}
  \begin{subfigure}[b]{0.32\textwidth}
    \begin{tikzpicture}[scale=0.49]
      \begin{axis}[ymode=log, grid style=dashed, xlabel={refinement}, ylabel={$\lVert u-u_h\rVert_{H^2(\Omega)}$}, title={$k=20$, $\alpha=10^{-2}$}, legend pos=south west, legend cell align={left}]
        \addplot table [x=refs, y=H2, col sep=comma] {assets/convergence_20.0_0.01_0.0.csv};
        \addplot[color=blue, mark=square] table [x=refs, y=H2, col sep=comma,] {assets/convergence_HCT_3_20.0_0.01_0.0.csv};
        \addplot [mark=*, color=red] table [x=refs, y=H2, col sep=comma] {assets/convergence_20.0_0.01_0.005.csv};
        \addplot[color=red, mark=square] table [x=refs, y=H2, col sep=comma,] {assets/convergence_HCT_3_20.0_0.01_0.005.csv};
        \addplot[gray,dashed,domain=1:4] {0.06*((1/16)^x)};
        \addplot[gray,dashed,domain=1:4] {4*((1/4)^x)};
        \legend{ARG $\beta=0$, HCT $\beta=0$, ARG $\beta=5\cdot 10^{-3}$, HCT $\beta=5\cdot 10^{-3}$}
      \end{axis}
    \end{tikzpicture}
  \end{subfigure}
  \begin{subfigure}[b]{0.32\textwidth}
    \begin{tikzpicture}[scale=0.48]
      \begin{axis}[ymode=log, grid style=dashed, xlabel={refinement}, ylabel={$\lVert u-u_h\rVert_{H^2(\Omega)}$}, title={$k=30$, $\alpha=10^{-2}$}, legend pos=south west, legend cell align={left}]
        \addplot table [x=refs, y=H2, col sep=comma] {assets/convergence_30.0_0.01_0.0.csv};
        \addplot[color=blue, mark=square] table [x=refs, y=H2, col sep=comma,] {assets/convergence_HCT_3_30.0_0.01_0.0.csv};
        \addplot [mark=*, color=red] table [x=refs, y=H2, col sep=comma] {assets/convergence_30.0_0.01_0.005.csv};
        \addplot[color=red, mark=square] table [x=refs, y=H2, col sep=comma,] {assets/convergence_HCT_3_30.0_0.01_0.005.csv};
        \addplot[gray,dashed,domain=1:4] {0.125*((1/16)^x)};
        \addplot[gray,dashed,domain=1:4] {8*((1/4)^x)};
        \legend{ARG $\beta=0$, HCT $\beta=0$, ARG $\beta=5\cdot 10^{-3}$, HCT $\beta=5\cdot 10^{-3}$}
      \end{axis}
    \end{tikzpicture}
  \end{subfigure}
  \\
  \hspace*{0.44cm}
  \begin{subfigure}[b]{0.3\textwidth}
    \includegraphics[width=0.9\textwidth]{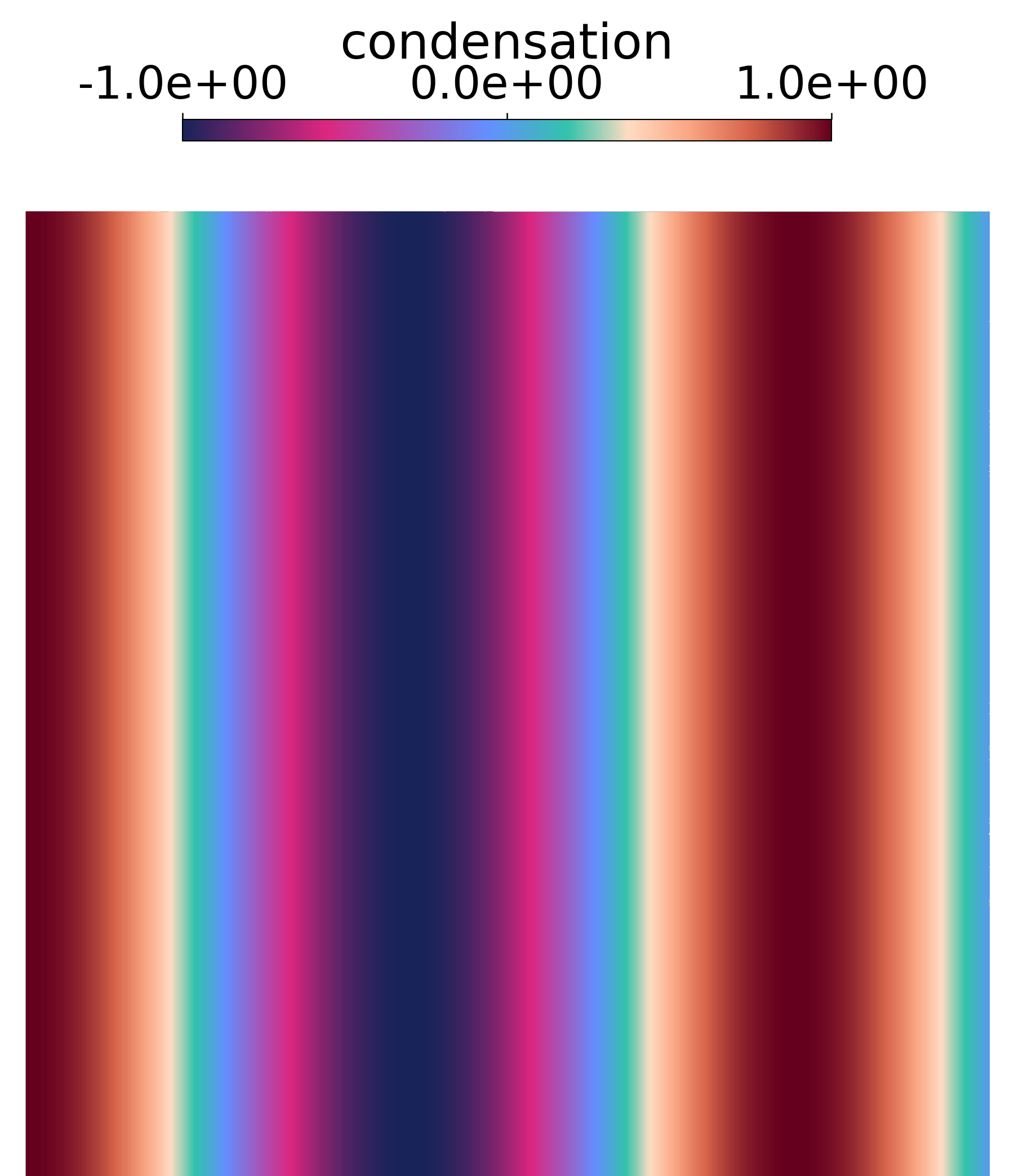}
  \end{subfigure}
  \hspace*{0.2cm}
  \begin{subfigure}[b]{0.3\textwidth}  
    \includegraphics[width=0.9\textwidth]{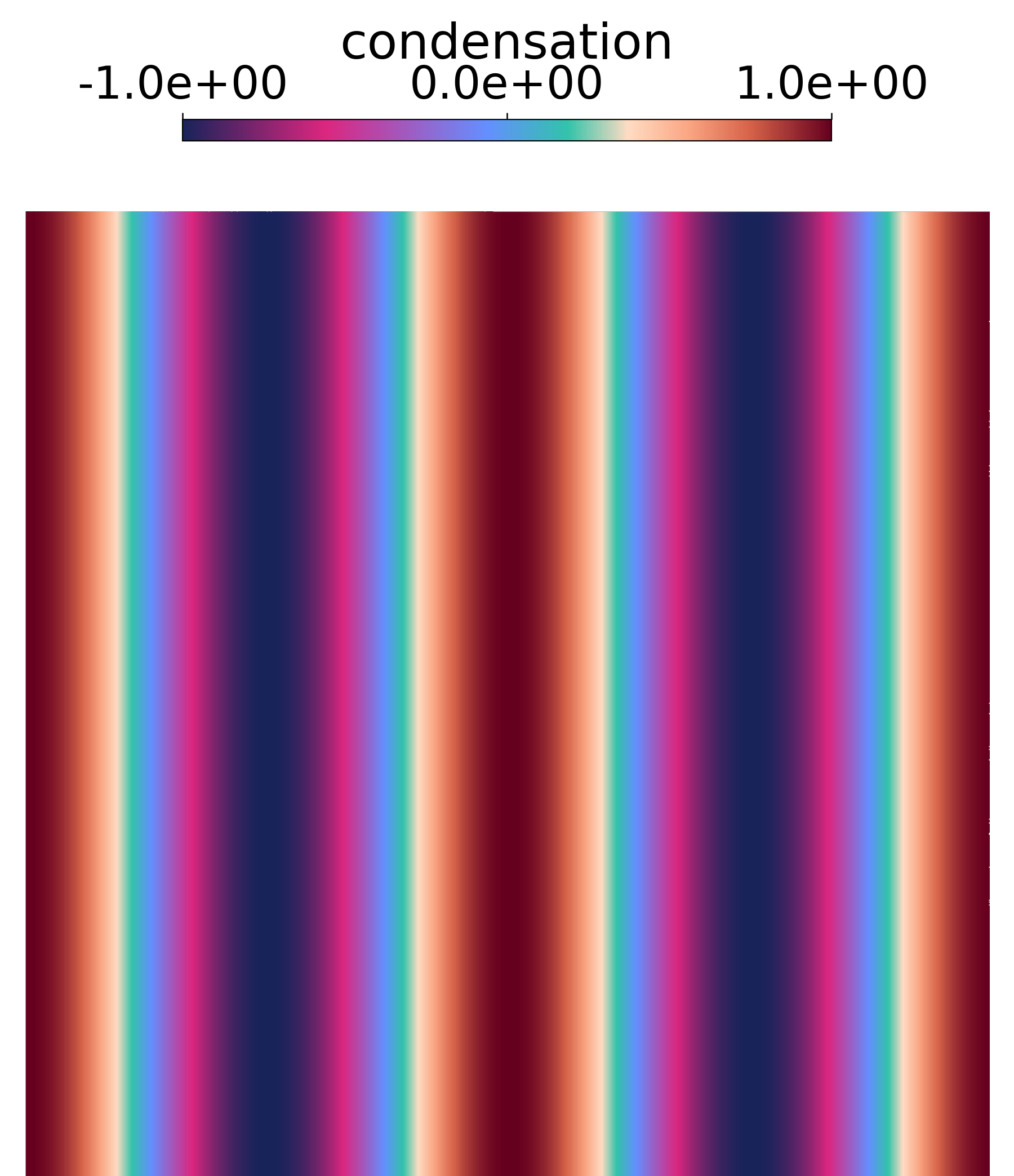}
  \end{subfigure}
  \hspace*{0.1cm}
  \begin{subfigure}[b]{0.3\textwidth}
    \includegraphics[width=0.9\textwidth]{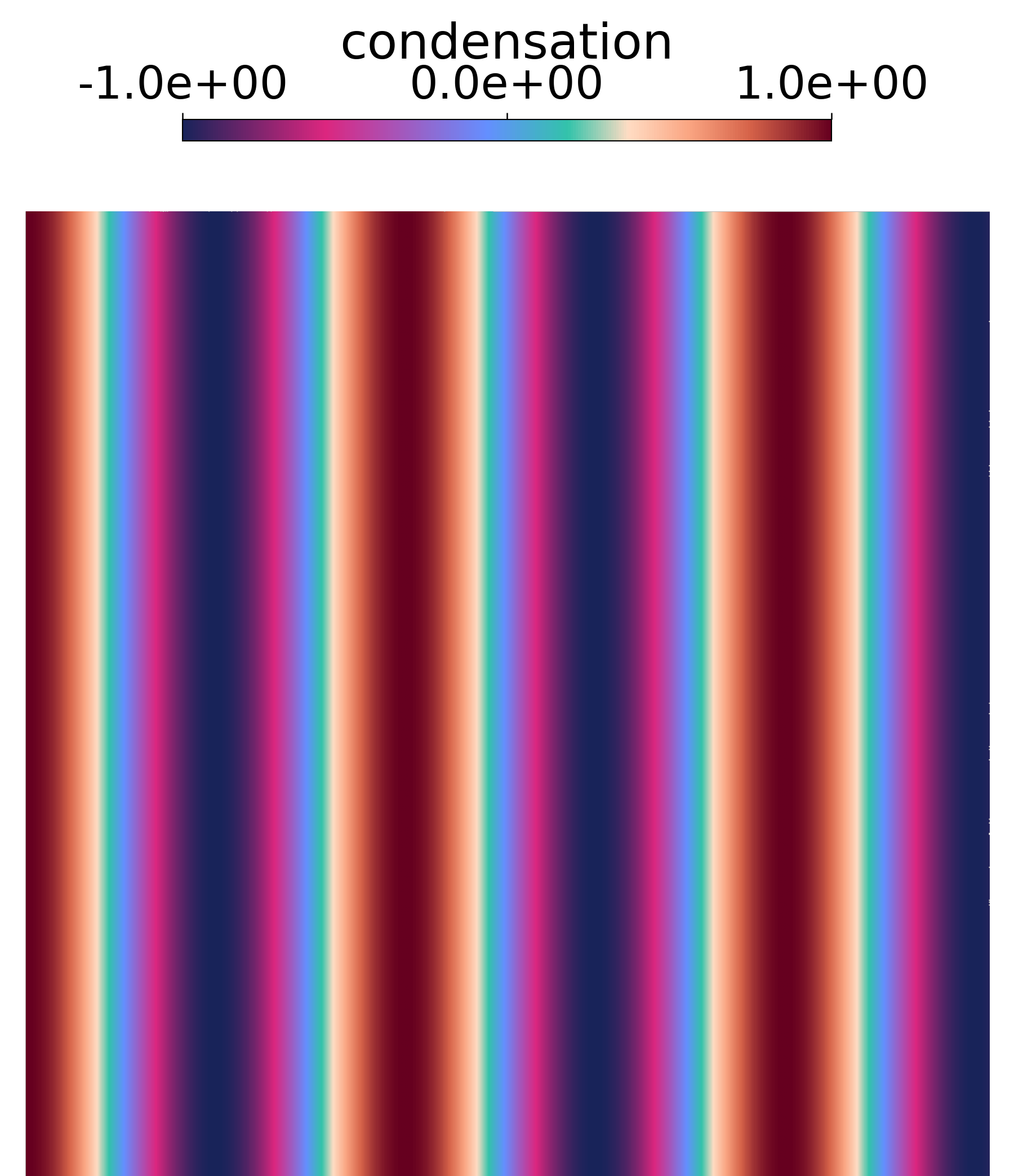}
  \end{subfigure}
  \hspace*{0.1cm}
  \centering
  \caption{The convergence of the $H^2$-norm of the error for the nematic Helmholtz--Korteweg equation for different values of $k$ (top row) and the corresponding manufactured solution (bottom row), where we display the real part of the condensation wave $u(\vec{x})$, i.e.\ $p(\vec{x},t) = \rho_0\left(1+\Re\{u(\vec{x}) e^{-i \omega t}\}\right)$, with $\rho(\vec{x},t)$ being the density at position $\vec{x}$ and time $t$ and $\rho_0$ the mean density of the fluid at rest.
} 
\label{fig:convergence}
\end{figure}

Plane waves are solutions of the Helmholtz equation that remain constant over a plane perpendicular to the direction of propagation. They are a natural choice for the numerical validation of the nematic Helmholtz--Korteweg equation, since they are also solutions of the latter, i.e.
\begin{equation}
  u(\vec{x}) = \exp\left(i(\vec{d}\cdot\vec{x})\right), \qquad \vec{d}\in \mathbb{C}.\label{eq:planewave}
\end{equation}
Substituting \eqref{eq:planewave} into \eqref{eq:nematic_helmholtz_korteweg}, it is easy to observe that plane waves are also exact solution of the nematic Helmholtz--Korteweg equation, under appropriate choice of the wave vector $\vec{d}$ \cite{farrell2024b}.
In particular, if $d \coloneqq \vert \vec{d} \vert$ satisfies the  following dispersion relation \eqref{eq:dispersion}, then the plane wave is also a solution of the nematic Helmholtz--Korteweg equation:
\begin{equation}
  \alpha d^4 + \beta d^2 (\vec{d}^\top\vec{n})^2+ d^2 - k^2  = 0.\label{eq:dispersion}
\end{equation}
In Figure \ref{fig:convergence} we present the convergence of the $H^2$-norm of the error for the nematic Helmholtz--Korteweg equation with boundary data constructed in such a way that the exact solution is a plane wave.
We observe that from Theorem \ref{thm:bestapprox} the error should converge with a rate only determined by the approximation properties of the finite element space.
In particular, it is well known, see for example \cite[Section 5.9]{BrennerScott}, that one can construct an interpolation operator from $H^5(\Omega)$ to the Argyris finite element space that converges with rate $h^3$ with respect to the $H^2$-norm \cite{Argyris}, i.e.
\begin{equation}
  \Vert u - I_h u \Vert_{H^2(\Omega)} \le C h^4 \Vert u \Vert_{H^5(\Omega)}.\label{eq:interpolation}
\end{equation}

Thus, combining \eqref{eq:interpolation} with Theorem \ref{thm:bestapprox} we expect the error to converge with rate $h^4$ with respect to the $H^2$-norm, as confirmed by Figure \ref{fig:convergence}.
Similarly, if we consider the Hsieh--Clough--Tocher macroelement \cite{HCT} we expect the error to converge with rate $h^2$ with respect to the $H^2$-norm, as confirmed by Figure \ref{fig:convergence}.
%\pftodo{Why does the $k=10$ figure have HCT results, but no others? And we don't discuss the HCT results at all. As far as I can work out this is the only place with HCT results, so please add HCT results to the other figures.}
\subsection{Anisotropic Gaussian pulse}
\begin{figure}[htbp]
  \begin{subfigure}[b]{0.24\textwidth}
    \centering
    \includegraphics[width=0.9\textwidth]{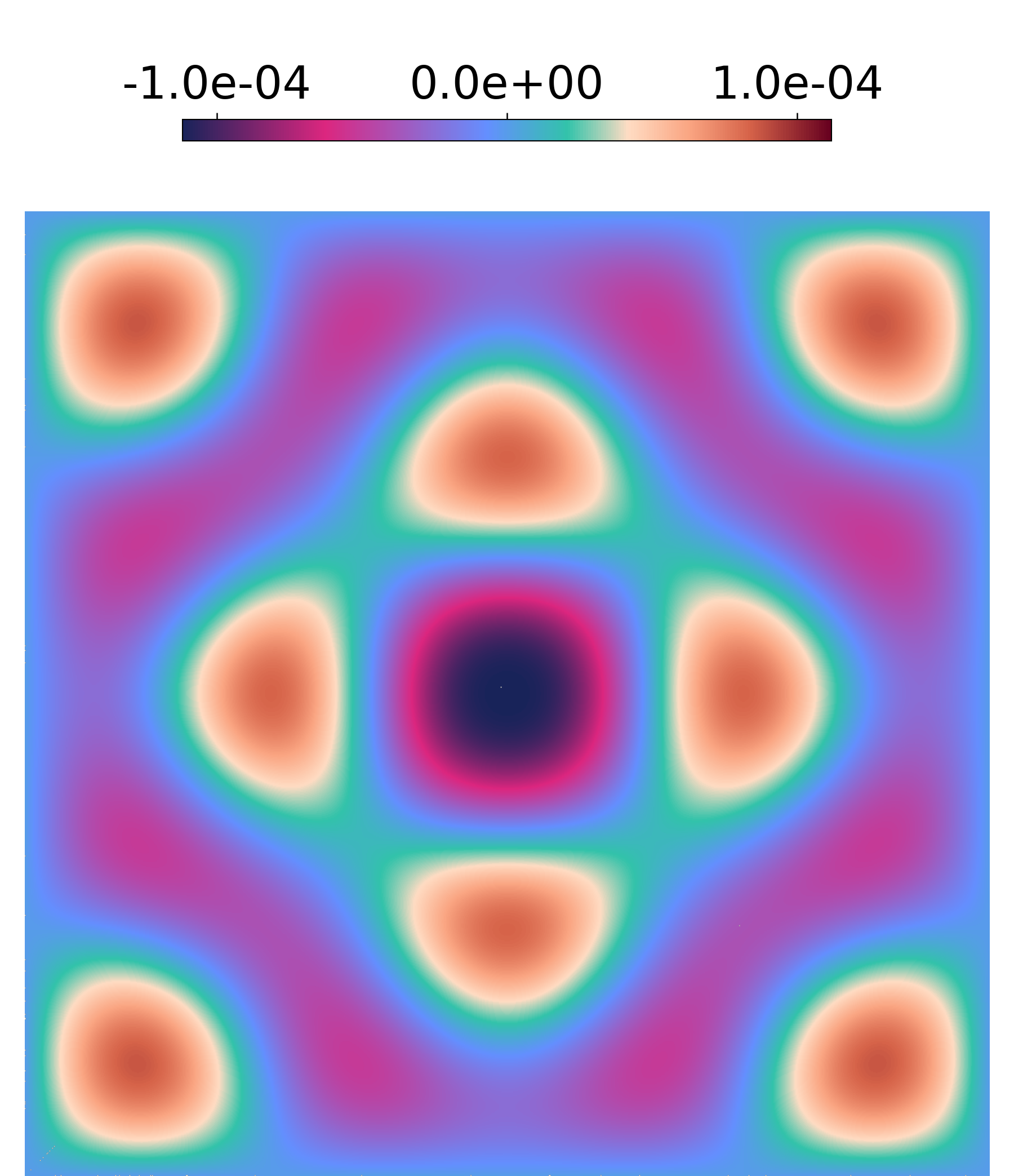}
    \subcaption{}
  \end{subfigure}
  \begin{subfigure}[b]{0.24\textwidth}
    \centering
    \includegraphics[width=0.9\textwidth]{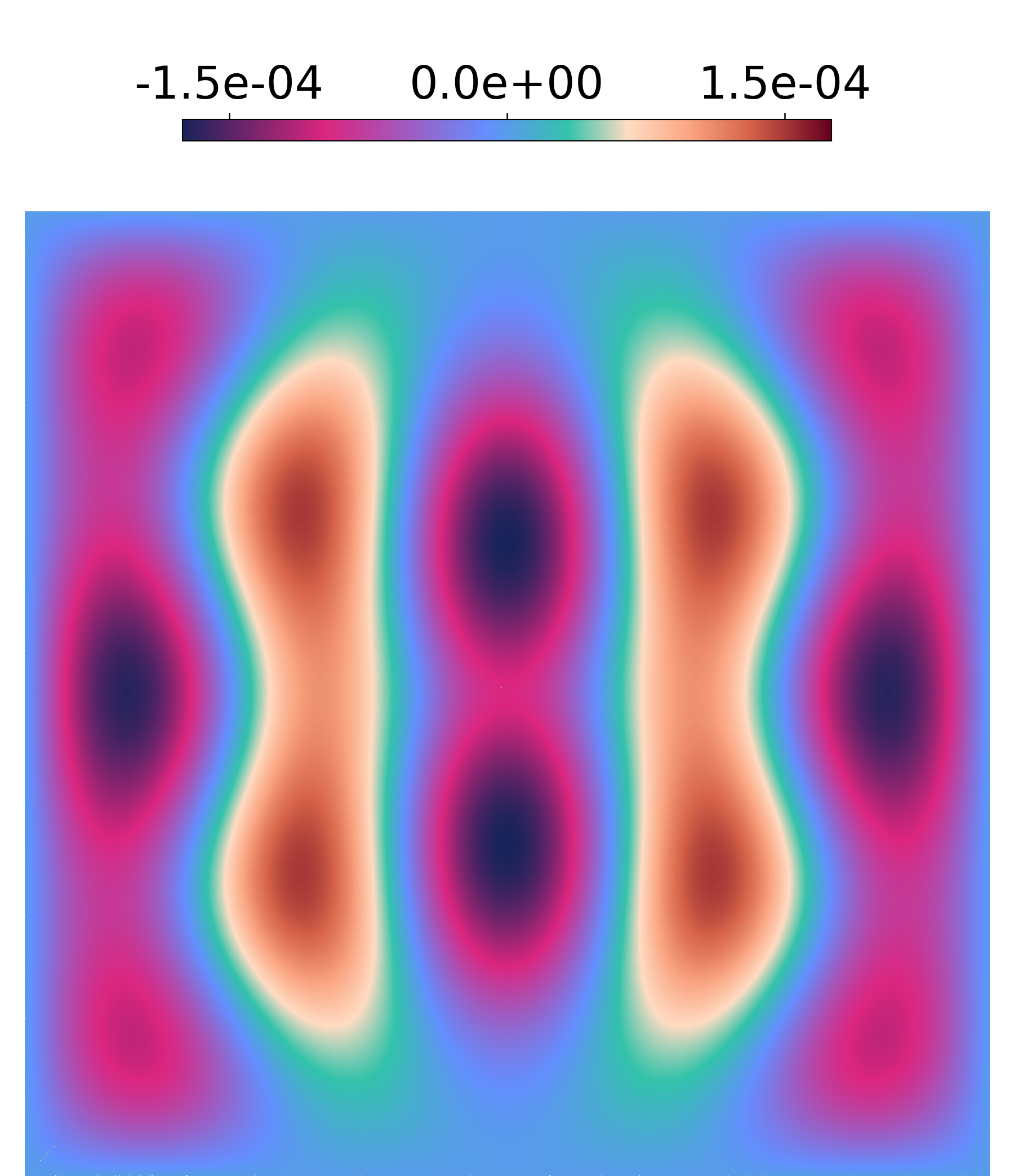}
    \subcaption{}
  \end{subfigure}
  \begin{subfigure}[b]{0.24\textwidth}
    \centering  
    \includegraphics[width=0.9\textwidth]{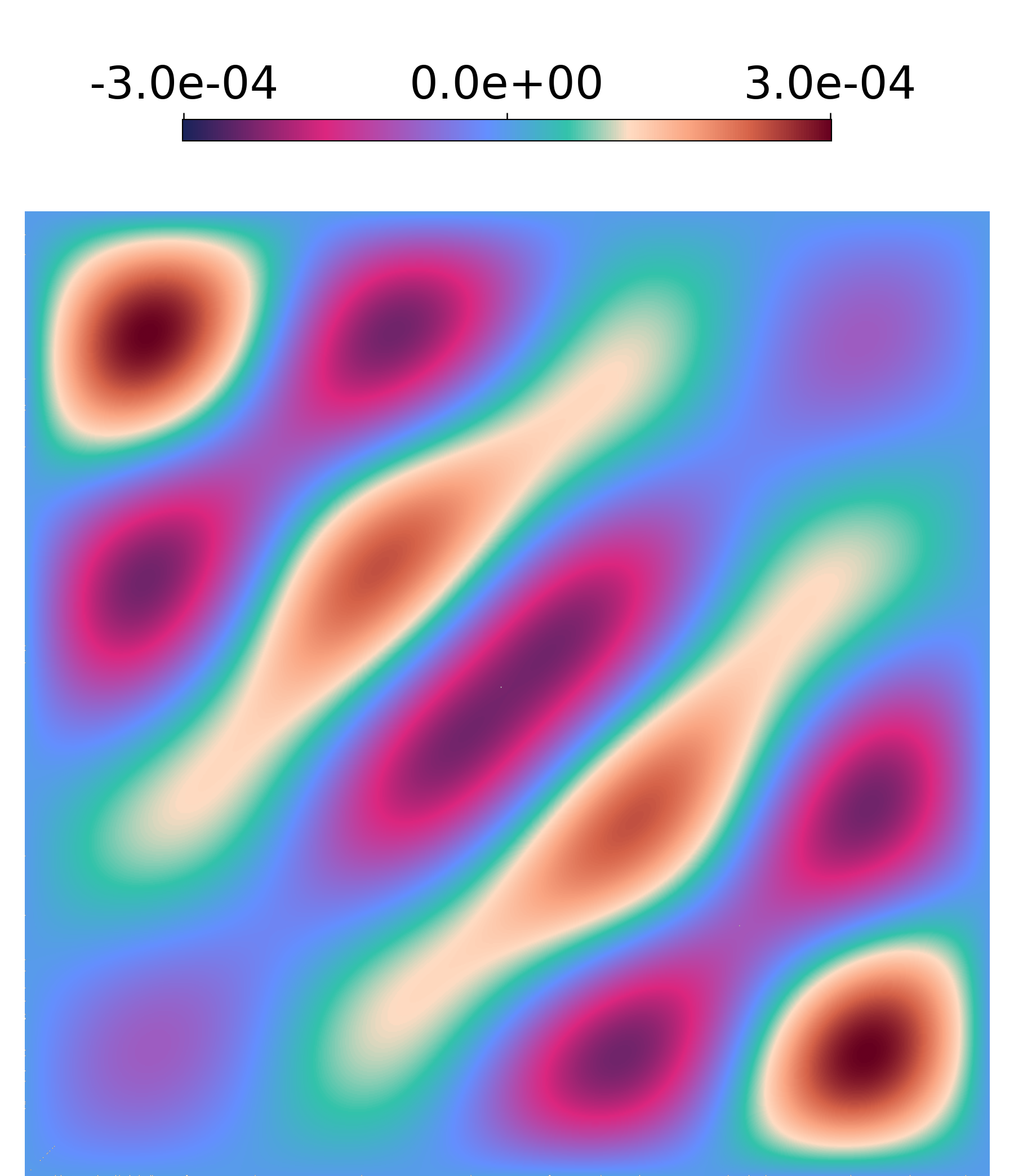}
    \subcaption{}
  \end{subfigure}
  \begin{subfigure}[b]{0.24\textwidth}
    \centering
    \includegraphics[width=0.9\textwidth]{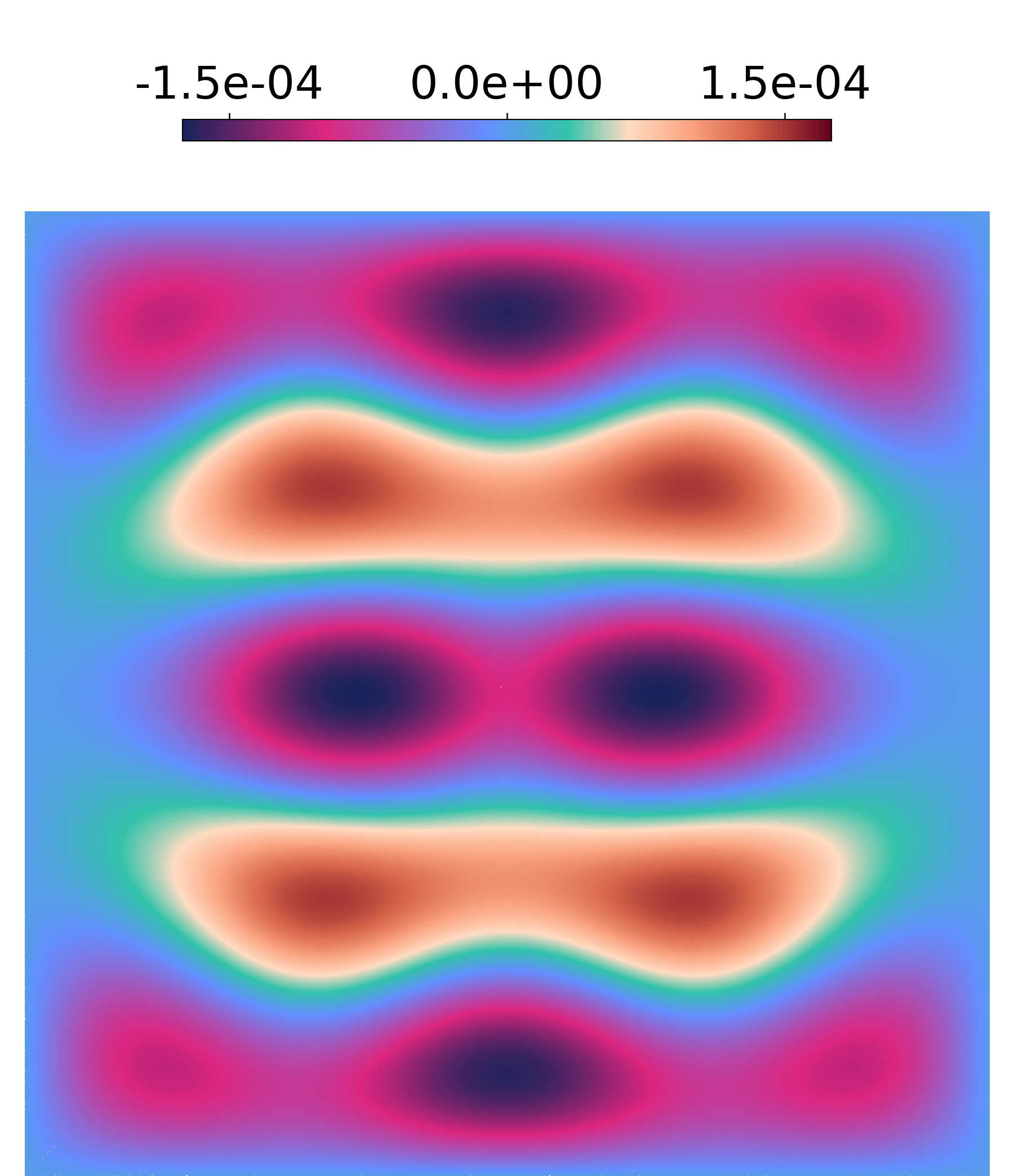}
    \subcaption{}
  \end{subfigure}
  \centering
\caption{The propagation of a symmetric Gaussian pulse by the nematic Helmholtz--Korteweg equation with sound soft boundary conditions when $\beta = 0$ (a) and when $\beta=5\cdot 10^{-3}$ and $\vec{n}$ is aligned with the $x$-axis (b), the diagonal (c) and the $y$-axis (d).
The parameters here are $k=50$, $\alpha=1\cdot 10^{-2}$.
}
\label{fig:GaussianPulseSS}
\end{figure}

\begin{figure}[htbp]
  \begin{subfigure}[b]{0.24\textwidth}
    \centering
    \includegraphics[width=0.9\textwidth]{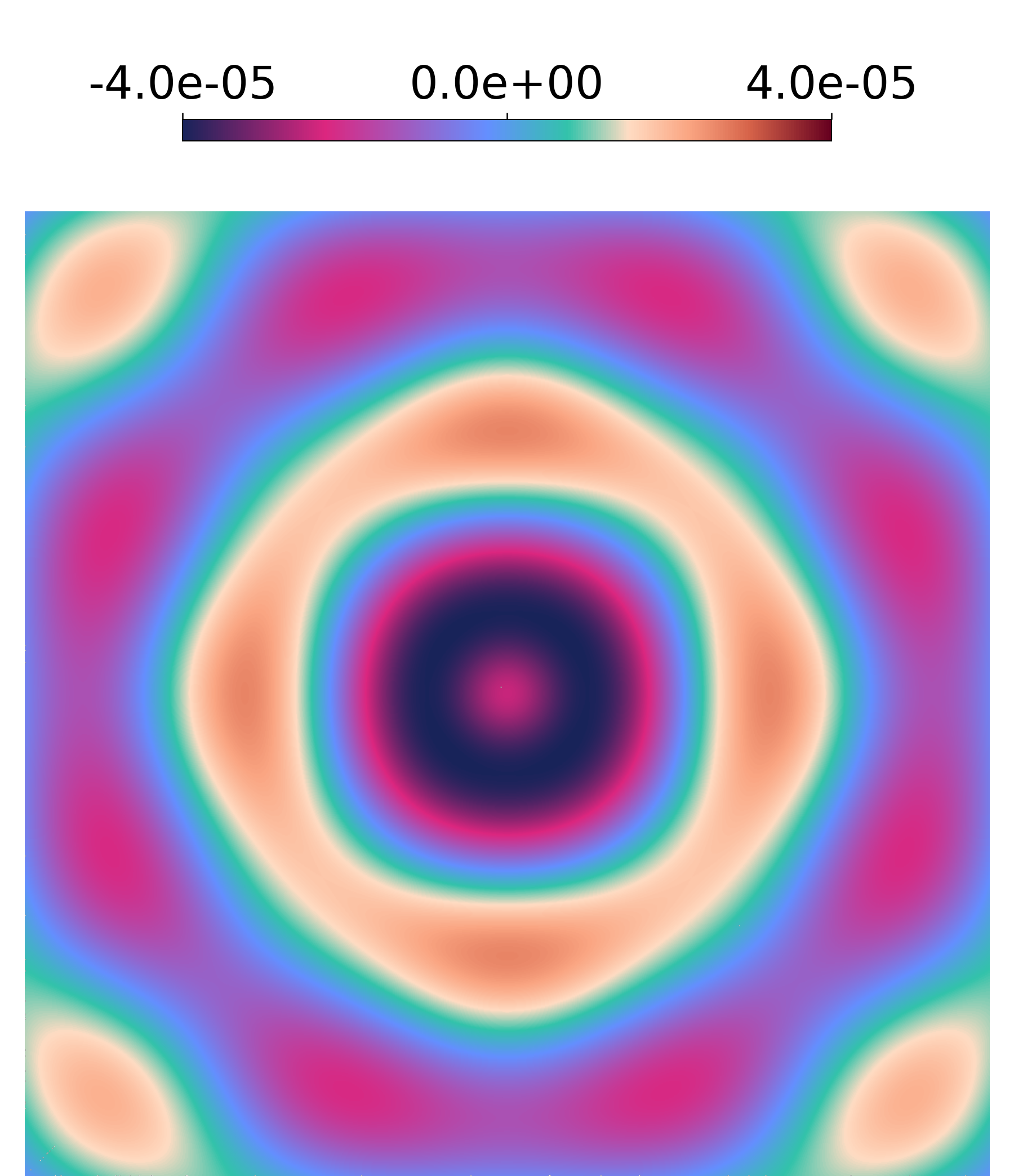}
    \subcaption{}
  \end{subfigure}
  \begin{subfigure}[b]{0.24\textwidth}
    \centering
    \includegraphics[width=0.9\textwidth]{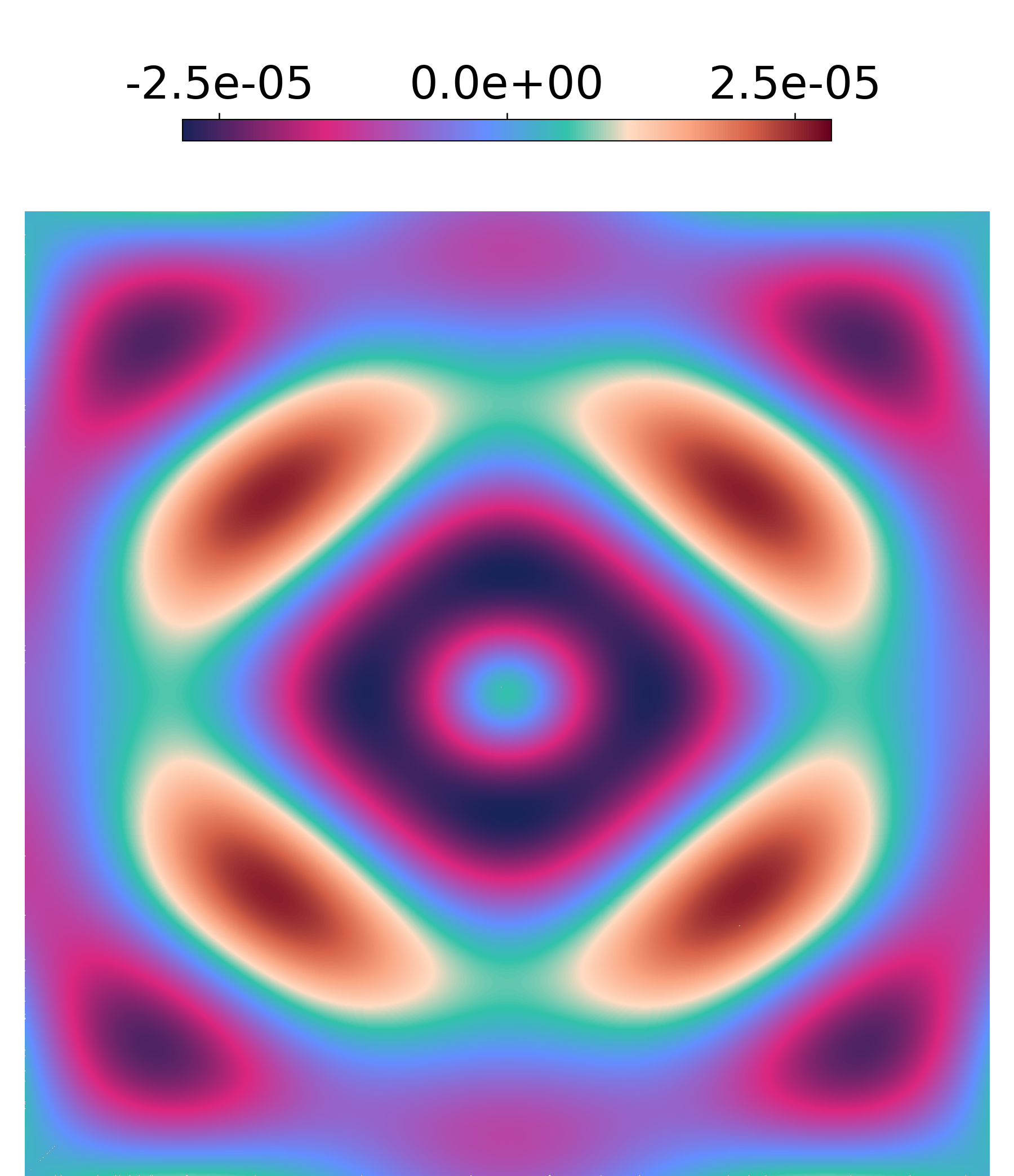}
    \subcaption{}
  \end{subfigure}
  \begin{subfigure}[b]{0.24\textwidth}
    \centering  
    \includegraphics[width=0.9\textwidth]{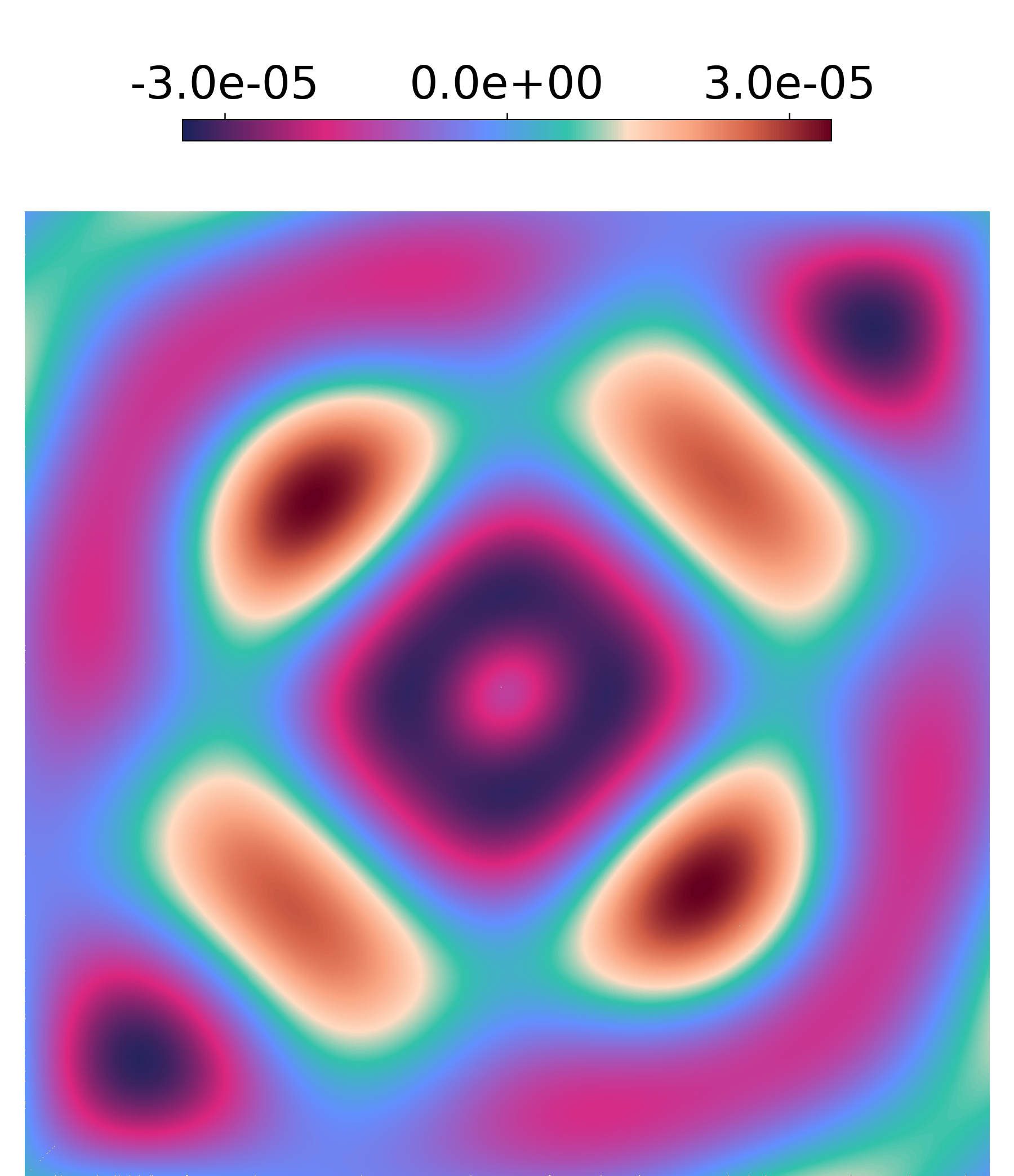}
    \subcaption{}
  \end{subfigure}
  \begin{subfigure}[b]{0.24\textwidth}
    \centering
    \includegraphics[width=0.9\textwidth]{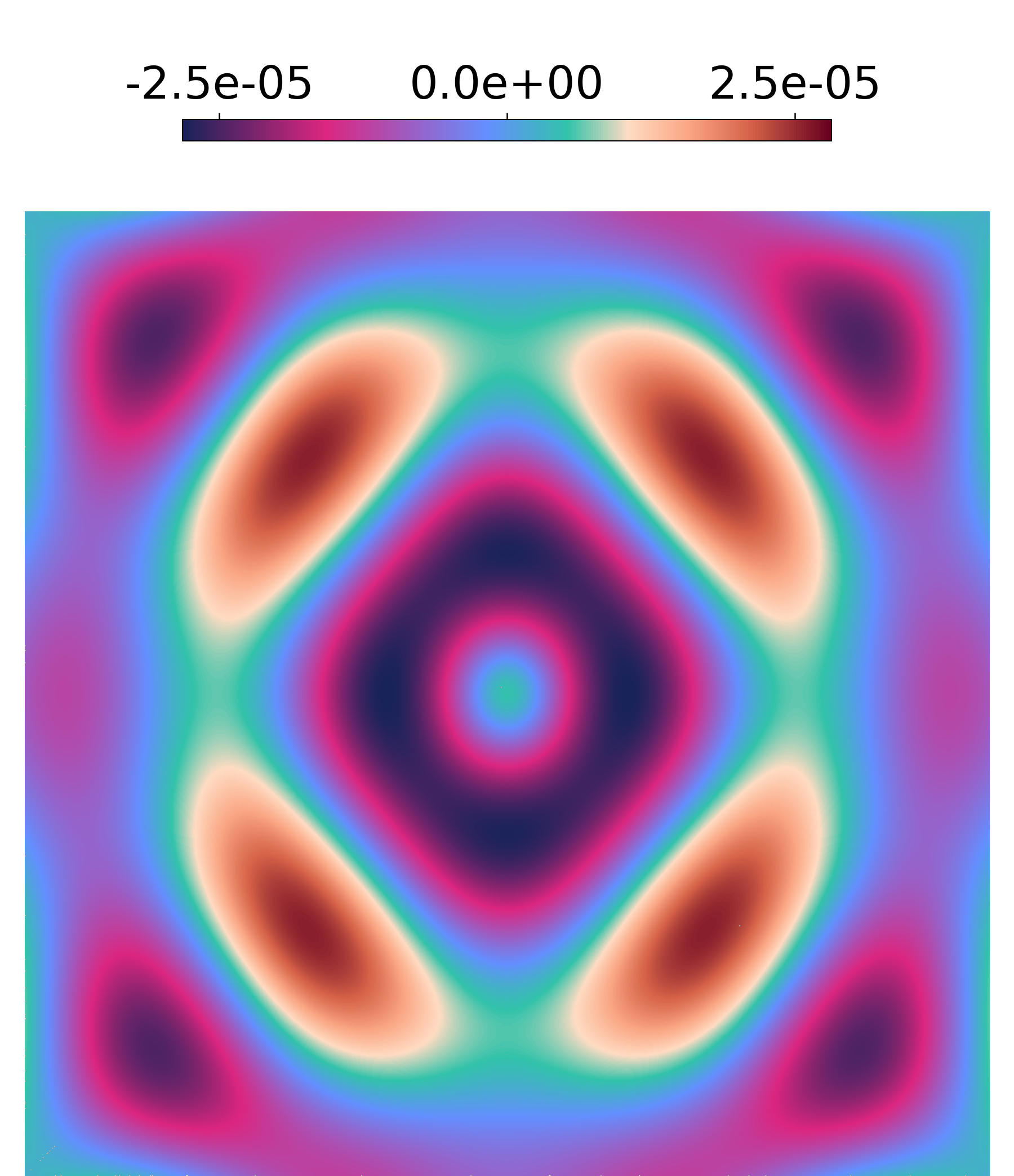}
    \subcaption{}
  \end{subfigure}
  \centering
\caption{The propagation of a symmetric Gaussian pulse by the nematic Helmholtz--Korteweg equation with impedance boundary conditions when $\beta = 0$ (a) and when $\beta=5\cdot 10^{-3}$ and $\vec{n}$ is aligned with the $x$-axis (b), the diagonal (c) and the $y$-axis (d).
The parameters here are $k=40$, $\alpha=1\cdot 10^{-2}$.
}
\label{fig:GaussianPulseImp}
\end{figure}

We study the anisotropic effect of the nematic term $\beta \nabla\cdot \nabla \left(\vec{n}^\top \mathcal{H} u \vec{n}\right)$ on the propagation of a symmetric Gaussian pulse of the form
\begin{equation}
  f(x,y) = \exp\left(-(40^2)\left[\left(x-\frac{1}{2}\right)^2 + \left(y-\frac{1}{2}\right)^2\right]\right).
\end{equation}
Several authors have studied the propagation of plane waves in the context of nematic liquid crystals \cite{farrell2024b,VirgaNematoacoustics,SelingerEtAll}.
In particular, it has been shown that the nematic term causes an anisotropic speed of propagation of plane waves, which is greater in the direction of the nematic director $\vec{n}$.
In Figure \ref{fig:GaussianPulseSS} we present numerical simulations of the propagation of a symmetric Gaussian pulse by the nematic Helmholtz--Korteweg equation with sound soft boundary conditions, and in Figure \ref{fig:GaussianPulseImp} we impose impedance boundary conditions.
Both figures show that the anisotropic speed of propagation observed for plane waves is also present for the Gaussian pulse, independent of the choice of boundary conditions.

\subsection{The Mullen--L\"{u}thi--Stephen experiment}
The Mullen--L\"{u}thi--Stephen experiment \cite{MullenLuti} consists of a planar acoustic wave propagating in a nematic liquid crystal, where the nematic director $\vec{n}$ kept fixed in the $x$-direction via a magnetic field.
The experiment reveals that the speed of propagation of the acoustic wave is anisotropic, being greater in the direction of the nematic director $\vec{n}$.
%\pef{Can we clarify how this result is different to the anisotropic propagation described in the previous subsection, i.e.~that we're no longer just considering plane waves?}
%\uztodo{Added a paragraph clarifying the difference between the two experiments.}

\begin{figure}
  \centering
  \includegraphics[width=0.37\textwidth]{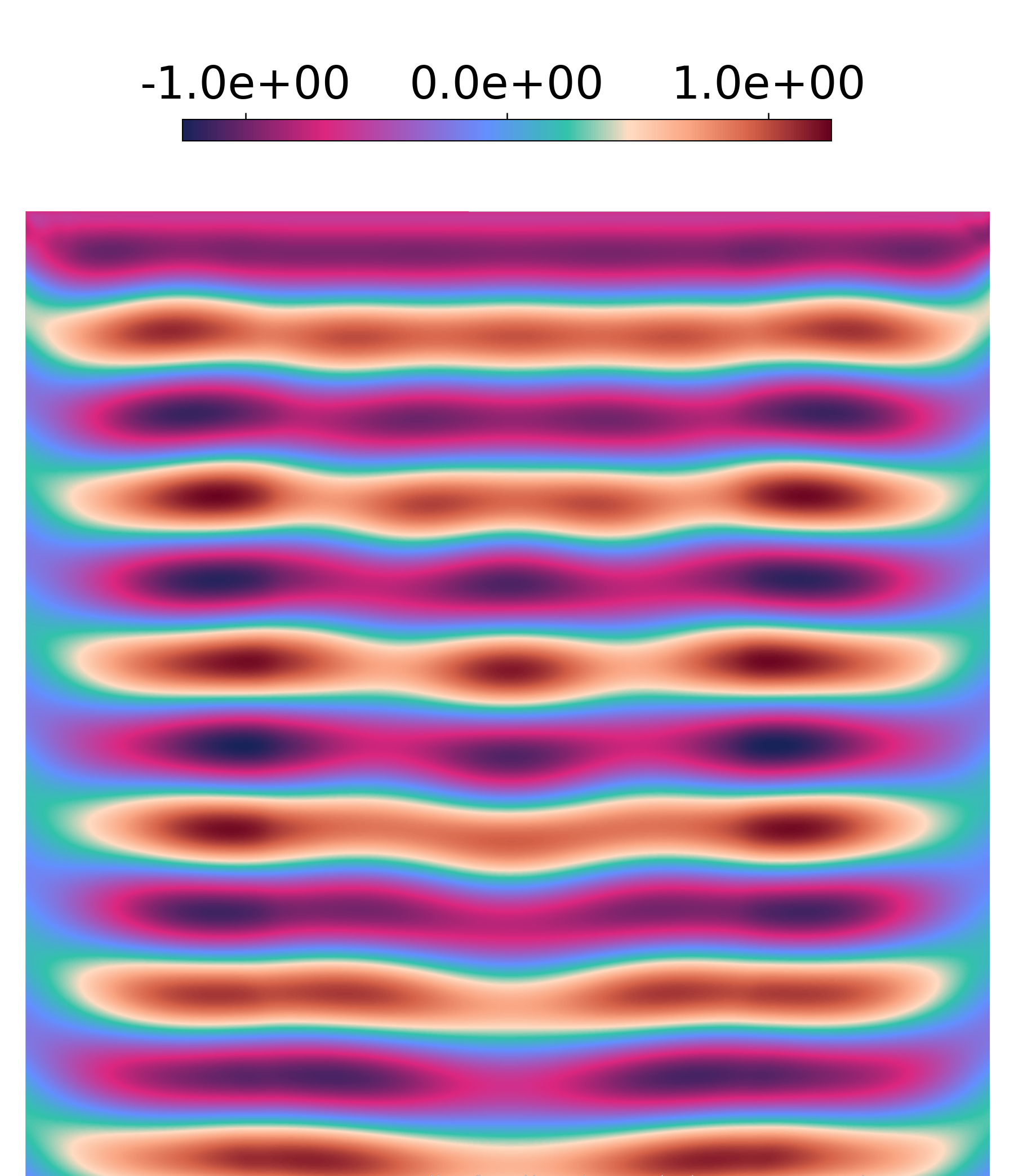}
  \hspace*{1.5em}
  \includegraphics[width=0.37\textwidth]{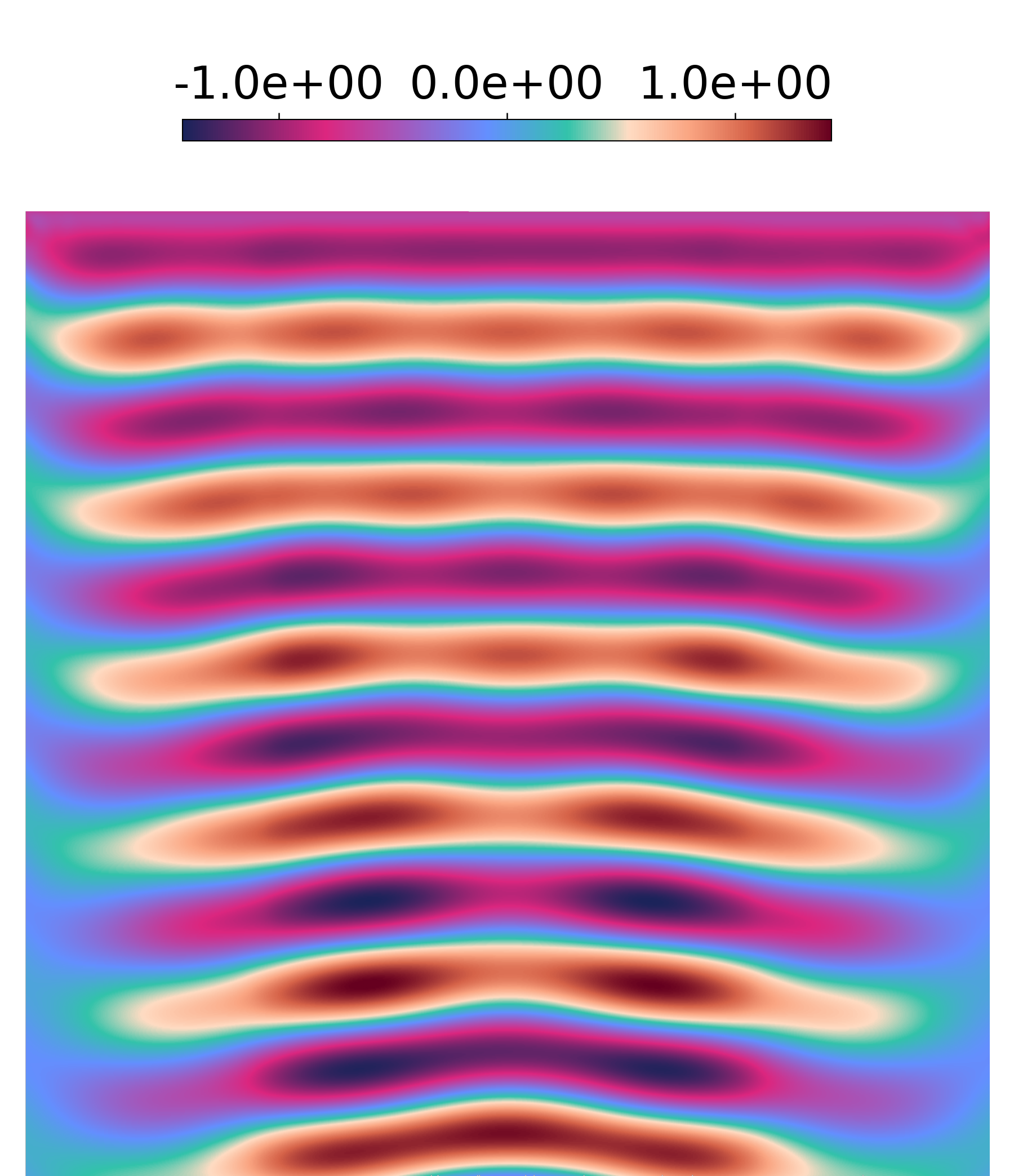}\\[1.5em]
  \includegraphics[width=0.32\textwidth]{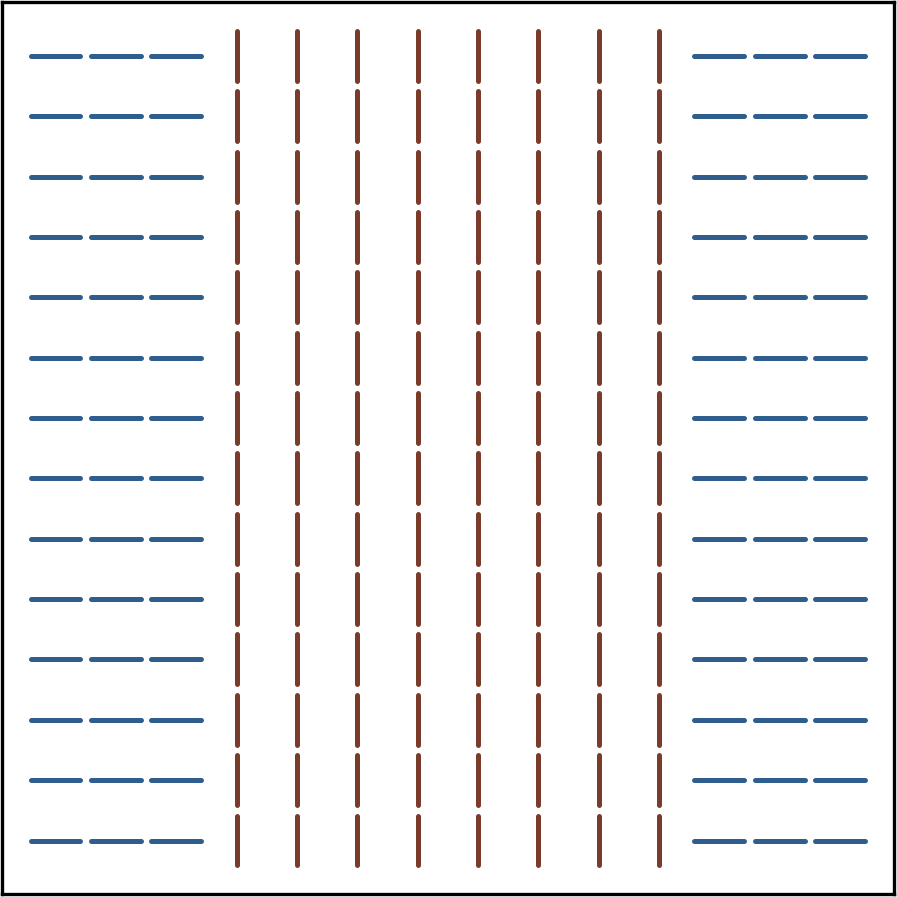}
  \hspace*{3.6em}
  \includegraphics[width=0.32\textwidth]{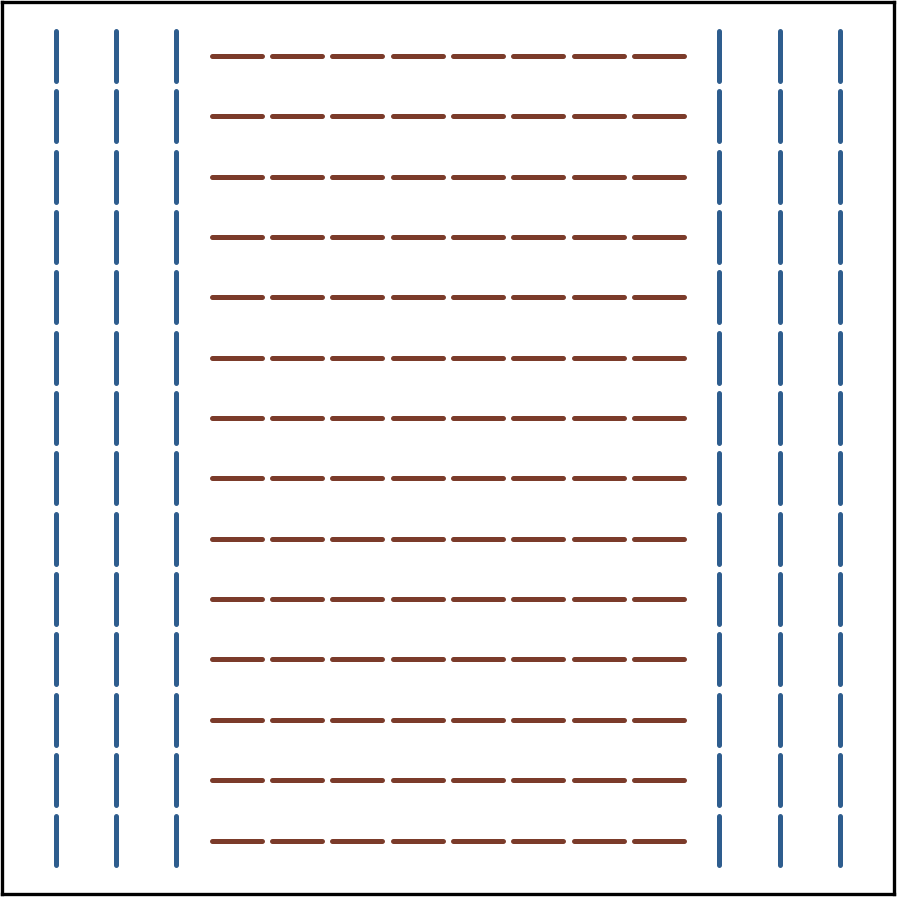}\\
  \caption{The anisotropic propagation of a planar acoustic wave in the modified Mullen--L\"{u}thi--Stephen experiment (top) and the corresponding nematic director $\vec{n}$ (bottom).
  The parameters are $k=40$, $\alpha=10^{-4}$ and $\beta= 5\cdot 10^{-5}$.}
  %\pftodo{I think we discussed the visualisation of this simulation in the PRE submission, and I have the same comments here: it is ugly.}
  \label{fig:MLS}
\end{figure}

We consider a variant of the Mullen--L\"{u}thi--Stephen experiment, where the nematic director $\vec{n}$ is kept fixed in the $x$-direction in the central region of the domain and in the $y$-direction in the outer region of the domain.
We then study the propagation of a planar acoustic wave coming from the top of the domain.
We can clearly see from Figure \ref{fig:MLS} that the speed of propagation of the acoustic wave is anisotropic, being greater in the central region of the domain where the nematic director $\vec{n}$ is aligned with the direction of propagation of the acoustic wave.
Furthermore, we can observe that the damping along the Robin boundary conditions imposed on the sides of the domain depends on the orientation of the nematic director $\vec{n}$.
In particular, the damping is greater in the region where the nematic director $\vec{n}$ is aligned with the direction of propagation of the acoustic wave.
This type of phenomena has been predicted from the analysis of the nematic Helmholtz--Korteweg equation \cite{farrell2024b}.

The difference between this experiment and the anisotropic propagation of plane waves discussed in the previous subsection is that rather than considering a Gaussian pulse propagating from the centre of the domain, here we consider a planar wave propagating from the top of the domain downwards,
and study how the orientation of the nematic director $\vec{n}$ affects the speed of propagation of the wave. In particular, while in the previous subsection we considered a nematic director $\vec{n}$ with constant orientation throughout the domain, here we consider a nematic director $\vec{n}$ with different orientations in different portions of the domain.

\begin{remark}
The discontinuity in the nematic director field presented in Figure \ref{fig:MLS} cannot be achieved if the nematic director field is a general solution of the Euler--Lagrange equations associated with the Oseen--Frank energy density, and should be regarded as an approximation of a solution of the Oseen--Frank model exhibiting a sharp transition between two constant orientations.
\end{remark}

\subsection{Tunable Resonators}\label{sec:TunableResonators}

Resonance occurs when the frequency of an acoustic wave coincides with one of the system's natural frequencies, corresponding to an eigenfunction of that system.
Under these conditions, the wave interferes constructively with itself, leading to a substantial amplification of the wave's amplitude.
This resonant behavior is not only a fundamental concept in wave physics but also a key principle in the design of acoustic resonators.
By carefully tuning geometrical and material parameters to match specific eigenfrequencies, engineers and scientists can enhance sound intensity, control wave propagation, and achieve highly efficient energy transfer within acoustic systems.

The nematic Helmholtz--Korteweg equation suggests that a nematic liquid crystal and other nematic materials can be used to design tunable acoustic resonators, since
the eigenvalues of the nematic Helmholtz--Korteweg equation depend on the orientation of the nematic director $\vec{n}$, as shown in \cite{farrell2024b}.
Thus, by changing the orientation of the nematic director $\vec{n}$, for example via an external electromagnetic field, it is possible to tune the eigenfrequencies of the system to be closer or further away from the frequency of the incoming acoustic wave, thus controlling the resonant behaviour of the system.

\begin{figure}[htbp]
  \centering
  \begin{subfigure}[b]{0.32\textwidth}
	  \centering
	  \includegraphics[width=\linewidth]{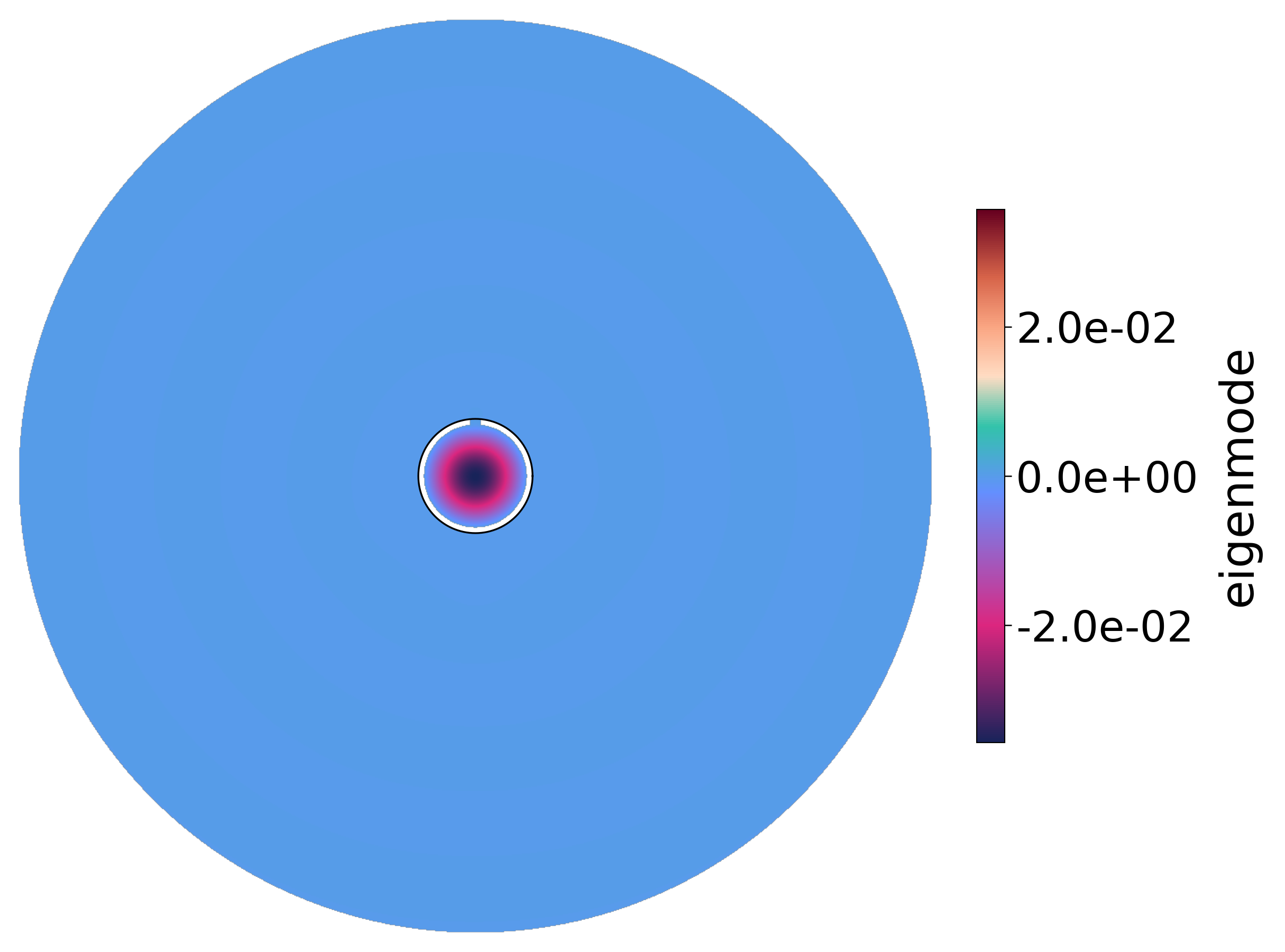}
	  \subcaption{}
  \end{subfigure}\hfill
  \begin{subfigure}[b]{0.32\textwidth}
	  \centering
	  \includegraphics[width=\linewidth]{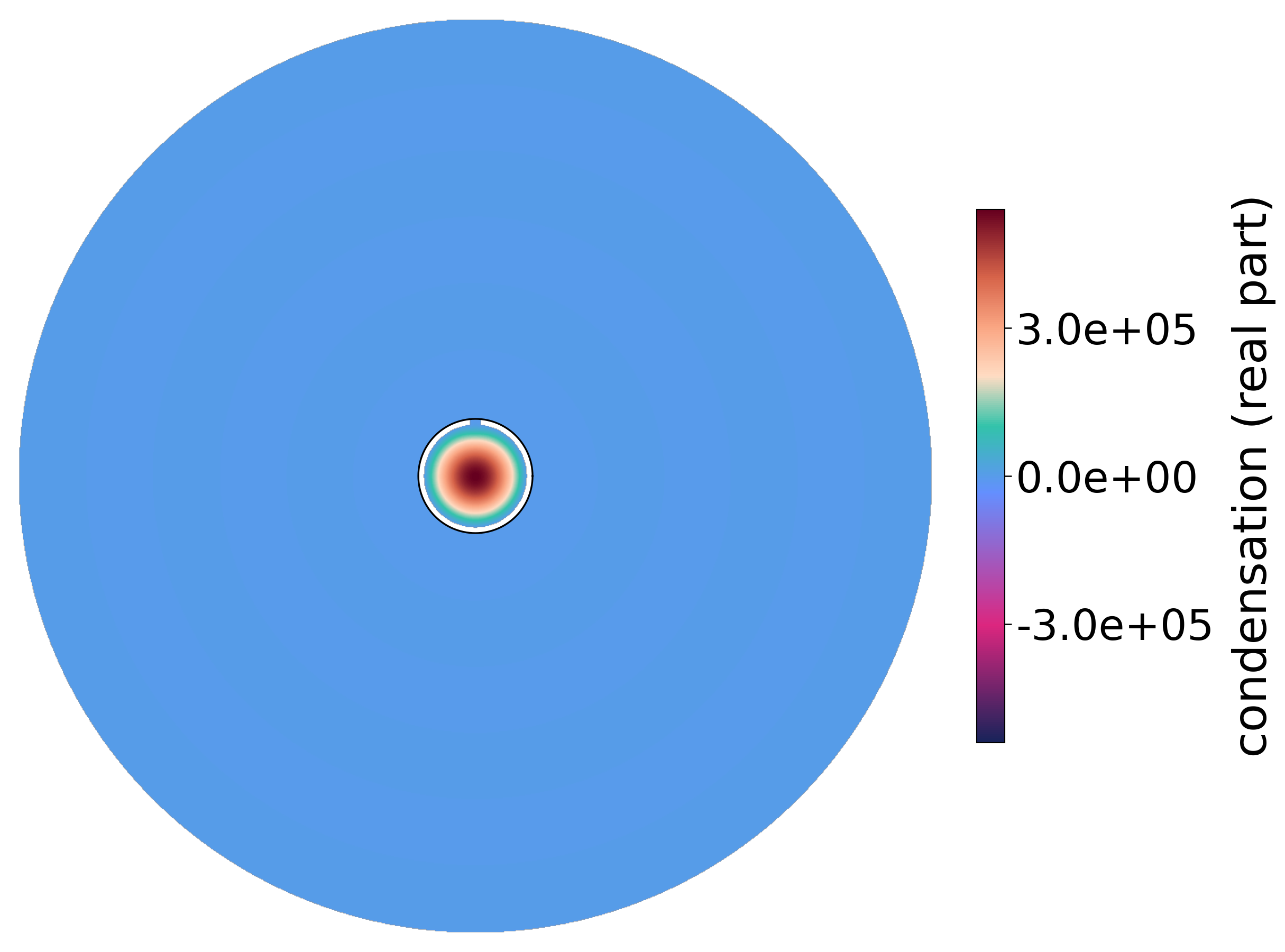}
	  \subcaption{}
  \end{subfigure}\hfill
  \begin{subfigure}[b]{0.32\textwidth}
	  \centering
	  \includegraphics[width=\linewidth]{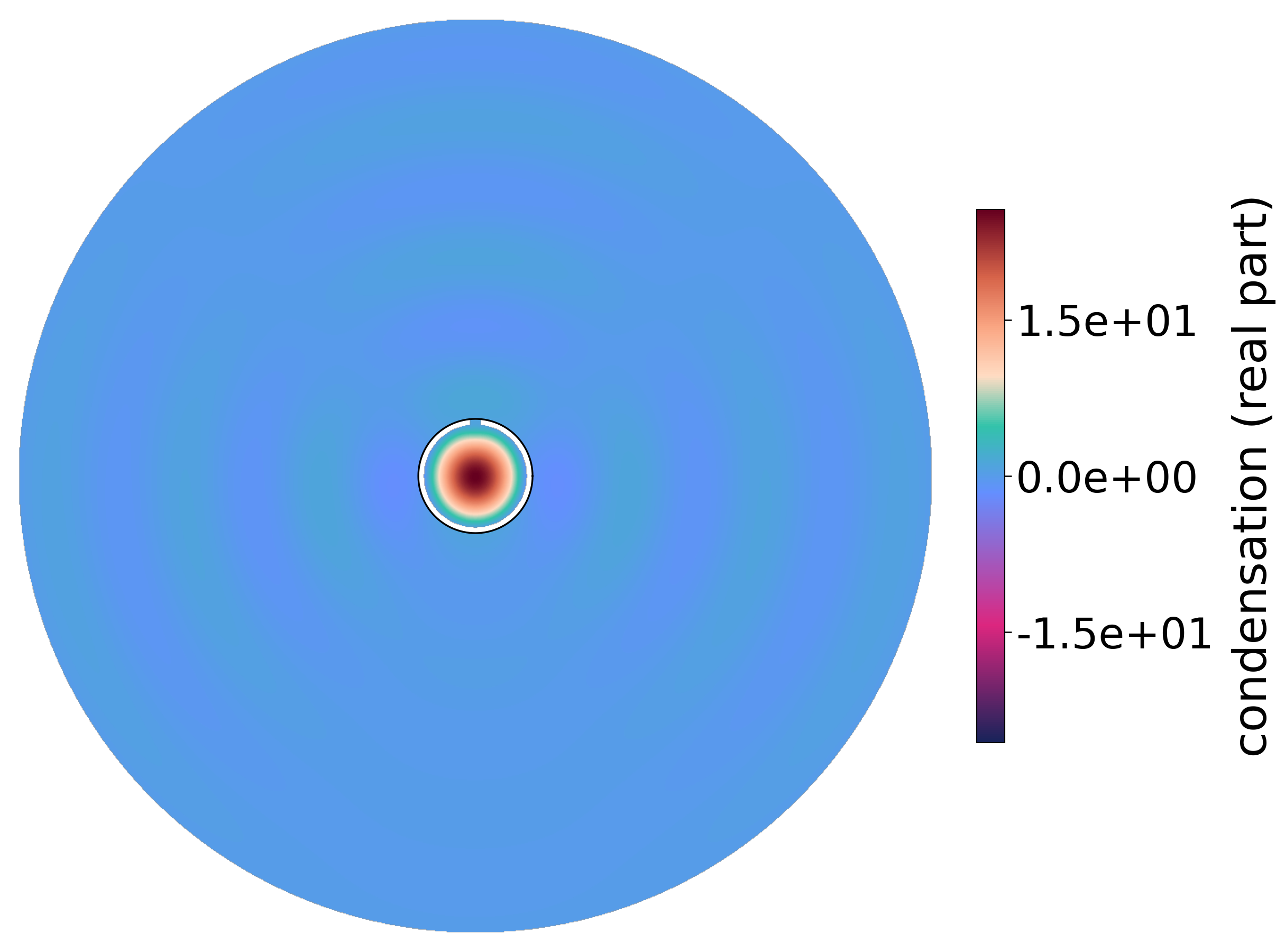}
	  \subcaption{}
  \end{subfigure}
  \caption{
	  (a) The eigenmode associated with the discrete eigenvalue \eqref{eq:numex:lambda} for $\vec{n} = (1,0)$.
	  (b) The scattered wave generated by incoming plane wave $\exp\left(i \vec{d}\cdot \vec{x}\right)$ with $\vec{d}=(0,\sqrt{Re(\lambda_h)})$ for a fixed nematic director $\vec{n} = (1,0)$
	  (c) The scattered wave generated by the same incoming plane wave for a fixed nematic director $\vec{n} = (0,1)$ respectively.
Notice the difference in scale: the middle figure has scale $10^4$ times larger.
  }
  \label{fig:Resonator}
\end{figure}

In Figure \ref{fig:Resonator} we present numerical simulations of a tunable acoustic resonator based on the nematic Helmholtz--Korteweg equation.
The simulation is obtained solving the following scattering problem for an incoming plane wave $u^{-}$:
\begin{equation}
  \begin{aligned}
      \label{eq:scatteringStrong}
    \alpha \Delta^2 u^+ +\nabla \cdot \nabla  \left(\vec{n}^\top (\mathcal{H}u^+) \vec{n}\right) - \Delta u^{+} -k^2 u^+ \! &=\! 0, &&\text{ in }\mathbb{R}^2\backslash\mathcal{R}, \\
      u^+ - u^-&= 0,&&\text{ on }\partial\mathcal{R} \\
      \alpha \Delta (u^+ - u^-) + \beta\vec{n}^\top \mathcal{H}(u^+ - u^-)\vec{n} &= 0, &&\text{ on } \partial\mathcal{R} \\
      \abs{\partial_{\abs{\vec{x}}}u^+-ik u^+}&=\mathcal{O}(\abs{\vec{x}}^{-\frac{1}{2}}), &&\abs{\vec{x}} \to \infty,
  \end{aligned}
\end{equation}
where $u^+$ is the scattered wave and $\mathcal{R}$ is the resonator's domain. We discretise the associated weak formulation with the Argyris finite element space. We truncate the unbounded domain $\mathbb{R}^2\backslash\mathcal{R}$ to a bounded domain via an adiabatic absorbing layer, as discussed in \cite{pml}. In particular, for $\vec{n} = (1,0)$ we set $u^{-}= \exp\left(i \vec{d}\cdot \vec{x}\right)$ with $\vec{d}=(0,\sqrt{Re(\lambda_h + \epsilon)})$, $\epsilon > 0$, where 
\begin{align}\label{eq:numex:lambda}
  \lambda_h = 21.829917701997985-2.4268903501127466\cdot 10^{-9}i
\end{align}
is a discrete eigenvalue of the full nematic Helmholtz--Korteweg operator $\mathcal{A}_h$ (i.e.~the discrete bilinear form $a_h(\cdot,\cdot)$ without the $-k^2(\cdot,\cdot)_{L^2}$ term, with the boundary conditions of the scattering configuration imposed via Nitsche on the scatterer boundary $\partial\mathcal{R}$ and a Robin condition on the truncation boundary) with $\alpha = 10^{-2}$ and $\beta = 5\cdot 10^{-3}$.
Note that we have to add a small perturbation $\epsilon$ to the eigenvalue to ensure that the assumption $k^2 \not \in \{ \lambda_h^{(i)} \}$ is satisfied.

From Figure \ref{fig:Resonator}(b) we can observe that if the nematic director $\vec{n}$ is aligned with the $x$-axis and the frequency of the incoming plane wave is close to the eigenfrequency of the system for the same nematic director orientation, then the system exhibits resonant behaviour, with a substantial amplification of the wave's amplitude. It is also possible to observe from Figure \ref{fig:Resonator}(c) that if the nematic director $\vec{n}$ is aligned with the $y$-axis, then the frequency of the incoming plane wave is no longer close to an eigenfrequency of the system, and the system does not exhibit resonant behaviour. This shows that by changing the orientation of the nematic director $\vec{n}$ it is possible to tune the resonant behaviour of the system, thus giving rise to a new class of tunable acoustic resonators.

%\subsection{Scattering by sound soft obstacles}
%In \cite{farrell2024b} the authors study the scattering of plane waves by sound soft circular obstacles in the context of the nematic Helmholtz--Korteweg equation.
%In the present section we present some numerical simulations of the same problem, where the scatterer has different shapes.
%To impose the radiating boundary conditions on the outer boundary of the domain, typical for scattering problems, we use an adiabatic absorbing layer, which we preferred over the perfectly matched layer due to its simplicity of implementation \cite{pml}.

\section*{Acknowledgements} The authors would like to thank Christoph Lehrenfeld and Ilaria Perugia for their valuable comments and suggestions.

\section*{Funding}
This work was funded by the Engineering and Physical Sciences Research Council
[grant numbers EP/R029423/1 and EP/W026163/1], the Donatio Universitatis Carolinae Chair ``Mathematical modelling of multicomponent systems'', and the UKRI Digital Research Infrastructure Programme through the Science and Technology Facilities Council's Computational Science Centre for Research Communities (CoSeC). The second author acknowledges funding through the DFG Grant 432680300 -- SFB 1456.

\appendix
\section{Traces}\label{sec:appendix:traces}
\tealchange{The boundary conditions in \eqref{eq:HKBC} and the boundary terms appearing in the sesquilinear form \eqref{eq:BFImpedance} involve the formal expressions
\begin{equation}\label{eq:biharmTraces}
  T_0 u \coloneqq \alpha \Delta u + \beta \vec{n}^\top (\mathcal{H}u) \vec{n}, \qquad
  T_1 u \coloneqq \alpha \nabla(\Delta u) \cdot \bm{\nu} + \beta \nabla\bigl(\vec{n}^\top (\mathcal{H}u) \vec{n}\bigr) \cdot \bm{\nu},
\end{equation}
which, for an arbitrary $u \in H^2(\Omega)$, do not admit a well-defined pointwise restriction to $\partial \Omega$. The purpose of this appendix is to show that, as is customary for the Neumann trace of the Laplace problem \cite[Lem.~2.1.3]{BBF13}, both $T_0 u$ and $T_1 u$ can be defined \emph{weakly}, that is, as dual objects whose action on test functions on $\partial \Omega$ is fixed by integration by parts in the interior.}

\tealchange{\paragraph{The Neumann trace by duality.}
For $\sigma \in H^1(\Omega)$, the Dirichlet trace $\operatorname{tr}(\sigma) \in H^{1/2}(\partial \Omega)$ is a bounded linear functional of $\sigma$ \cite{Leoni}. By contrast, for a generic $\sigma \in H^1(\Omega)$, the gradient $\nabla \sigma$ lies only in $L^2(\Omega;\mathbb{R}^d)$, and so the pointwise restriction $\nabla \sigma \cdot \bm{\nu} \vert_{\partial \Omega}$ has no intrinsic meaning. The standard remedy is to require the additional regularity $\nabla \sigma \in H(\text{div},\; \Omega)$, i.e.\ $\Delta \sigma \in L^2(\Omega)$, and then, following \cite[Lem.~2.1.3]{BBF13}, to \emph{define} the Neumann trace $\operatorname{tr}_{\bm{\nu}}(\sigma) \in H^{-1/2}(\partial \Omega)$ by pairing against test functions on the whole domain: $\operatorname{tr}_{\bm{\nu}}(\sigma)$ is the unique element such that
\begin{equation}\label{eq:NeumannTraceDuality}
  \langle \operatorname{tr}_{\bm{\nu}}(\sigma),\, \operatorname{tr}(v) \rangle_{H^{-1/2}(\partial \Omega),\,H^{1/2}(\partial \Omega)}
  \coloneqq \int_{\Omega} \nabla \sigma \cdot \nabla v \dx + \int_{\Omega} \Delta \sigma \, v \dx
  \qquad \text{for all } v \in H^1(\Omega).
\end{equation}
The right-hand side of \eqref{eq:NeumannTraceDuality} depends on $v$ only through $\operatorname{tr}(v) \in H^{1/2}(\partial \Omega)$ (the difference of two $v$'s with the same trace lies in $H^1_0(\Omega)$, on which the right-hand side vanishes by Green's identity), and by the surjectivity of $\operatorname{tr} : H^1(\Omega) \to H^{1/2}(\partial \Omega)$ this fixes $\operatorname{tr}_{\bm{\nu}}(\sigma)$ uniquely. It yields a bounded operator $\operatorname{tr}_{\bm{\nu}} : \{ \sigma \in H^1(\Omega) : \Delta \sigma \in L^2(\Omega) \} \to H^{-1/2}(\partial \Omega)$ \cite[Lem.~2.1.3]{BBF13}. The Neumann trace is thus not the restriction of a pointwise object but a \emph{dual} object.
}

\tealchange{\paragraph{The biharmonic traces $T_0 u$ and $T_1 u$ by duality.}
We now apply the same principle to \eqref{eq:biharmTraces}. The key observation is that, although for an arbitrary $u \in H^2(\Omega)$ neither $T_0 u\vert_{\partial \Omega}$ nor $T_1 u\vert_{\partial \Omega}$ make sense, any $u \in H^2(\Omega)$ solving \eqref{eq:nematic_helmholtz_korteweg} distributionally enjoys the additional interior regularity
\begin{equation}\label{eq:LapT0Reg}
  \Delta T_0 u
  = \alpha \Delta^2 u + \beta \nabla \cdot \nabla\bigl(\vec{n}^\top (\mathcal{H}u) \vec{n}\bigr)
  = f + \Delta u + k^2 u \in L^2(\Omega),
\end{equation}
because $f \in L^2(\Omega)$ by assumption, which is strictly more than the natural dual regularity $H^{-2}(\Omega)$. Consequently $T_0 u \in L^2(\Omega)$ with $\Delta T_0 u \in L^2(\Omega)$, i.e.\ $T_0 u \in H(\Delta,\Omega) \coloneqq \{ w \in L^2(\Omega) : \Delta w \in L^2(\Omega)\}$. We emphasise that we do \emph{not} require $\nabla T_0 u \in H(\operatorname{div},\Omega)$: such a condition would amount to $\nabla T_0 u \in L^2(\Omega;\mathbb{R}^d)$, a regularity not granted by $u \in H^2(\Omega)$ alone---the interior PDE only ensures that the divergence $\Delta T_0 u$ lies in $L^2(\Omega)$. The weaker regularity $T_0 u \in H(\Delta,\Omega)$ suffices because the construction below relies on Green's \emph{second} identity, placing both Laplacians on the smooth side of each duality pairing: for any sufficiently smooth $w$,
\begin{equation}\label{eq:biharmGreen}
  \int_{\Omega} (\Delta T_0 u)\, w \dx
  = \int_{\Omega} T_0 u \, \Delta w \dx
  \;+\; \int_{\partial \Omega} T_1 u \, w \ds
  \;-\; \int_{\partial \Omega} T_0 u \, (\nabla w \cdot \bm{\nu}) \ds.
\end{equation}
Identity \eqref{eq:biharmGreen} plays for our boundary operators the same role that the divergence theorem plays in \eqref{eq:NeumannTraceDuality}, and motivates the following two duality definitions.}

\tealchange{\emph{The trace $T_0 u \in H^{-1/2}(\partial \Omega)$.} We pair $T_0 u$ against the Neumann trace of an arbitrary test function $v \in H^2(\Omega)$ with $\operatorname{tr}(v) = 0$: the trace $T_0 u$ is the unique element of $H^{-1/2}(\partial \Omega)$ such that
\begin{equation}\label{eq:T0traceDef}
  \langle T_0 u,\, \operatorname{tr}_{\bm{\nu}}(v) \rangle_{H^{-1/2}(\partial \Omega),\, H^{1/2}(\partial \Omega)}
  \coloneqq \int_{\Omega} T_0 u \, \Delta v \dx - \int_{\Omega} (\Delta T_0 u)\, v \dx
  \qquad \text{for all } v \in H^2(\Omega) \text{ with } \operatorname{tr}(v) = 0.
\end{equation}
The right-hand side of \eqref{eq:T0traceDef} depends on $v$ only through $\operatorname{tr}_{\bm{\nu}}(v) \in H^{1/2}(\partial \Omega)$: by Green's identity \eqref{eq:biharmGreen} applied with $\operatorname{tr}(v) = 0$, the difference between two admissible $v$'s with the same $\operatorname{tr}_{\bm{\nu}}$ contributes zero. The surjectivity of $\operatorname{tr}_{\bm{\nu}}$ from $\{ v \in H^2(\Omega) : \operatorname{tr}(v) = 0\}$ onto $H^{1/2}(\partial \Omega)$ then fixes $T_0 u$ uniquely.}

\tealchange{\emph{The trace $T_1 u \in H^{-3/2}(\partial \Omega)$.} We pair $T_1 u$ against the Dirichlet trace of an arbitrary test function $v \in H^2(\Omega)$ with $\operatorname{tr}_{\bm{\nu}}(v) = 0$: the trace $T_1 u$ is the unique element of $H^{-3/2}(\partial \Omega)$ such that
\begin{equation}\label{eq:T1traceDef}
  \langle T_1 u,\, \operatorname{tr}(v) \rangle_{H^{-3/2}(\partial \Omega),\, H^{3/2}(\partial \Omega)}
  \coloneqq \int_{\Omega} (\Delta T_0 u)\, v \dx - \int_{\Omega} T_0 u \, \Delta v \dx
  \qquad \text{for all } v \in H^2(\Omega) \text{ with } \operatorname{tr}_{\bm{\nu}}(v) = 0.
\end{equation}
The right-hand side of \eqref{eq:T1traceDef} depends on $v$ only through $\operatorname{tr}(v) \in H^{3/2}(\partial \Omega)$, again by Green's identity \eqref{eq:biharmGreen}, and the surjectivity of $\operatorname{tr}$ from $\{ v \in H^2(\Omega) : \operatorname{tr}_{\bm{\nu}}(v) = 0\}$ onto $H^{3/2}(\partial \Omega)$ fixes $T_1 u$ uniquely.}

\tealchange{By \eqref{eq:LapT0Reg} and the Cauchy--Schwarz inequality, the right-hand sides of \eqref{eq:T0traceDef}--\eqref{eq:T1traceDef} depend continuously on $u$, so these definitions yield bounded linear operators
\begin{equation*}
  T_0 : \{ u \in H^2(\Omega) : \Delta T_0 u \in L^2(\Omega)\} \to H^{-1/2}(\partial \Omega), \qquad
  T_1 : \{ u \in H^2(\Omega) : \Delta T_0 u \in L^2(\Omega)\} \to H^{-3/2}(\partial \Omega),
\end{equation*}
in exact analogy with the Neumann trace \eqref{eq:NeumannTraceDuality}. We emphasise that, unlike the standard Dirichlet trace, neither of these biharmonic traces is compact: they are genuinely dual objects whose action is fixed only through the integration-by-parts identity \eqref{eq:biharmGreen}, and not through any pointwise restriction.
}

\tealchange{\begin{remark}
We emphasise that there is neither a gain in regularity for the solution $u$ of \eqref{eq:HKBC} nor any requirement on the regularity of $\partial \Omega$ beyond what is needed for the standard trace and inverse trace theorems on $H^1(\Omega)$ and $H^2(\Omega)$. We only seek $u \in H^2(\Omega)$ and use the equation in the interior of $\Omega$ to obtain the additional regularity \eqref{eq:LapT0Reg} that enables the duality definitions \eqref{eq:T0traceDef}--\eqref{eq:T1traceDef}. A separate elliptic regularity argument can be used to upgrade $u$ to $H^4(\Omega)$ under appropriate assumptions on $f$ and $\partial \Omega$ \cite{Grisvard}.
\end{remark}
}

\tealchange{Finally, we note that we could not have looked for a solution $u \in H^2(\Omega)$ directly in the subspace $\{ v \in H^2(\Omega) : \alpha \Delta^2 v + \beta \nabla \cdot \nabla (\vec{n}^\top (\mathcal{H}v) \vec{n}) \in L^2(\Omega)\}$: while this set is a closed subspace of $H^2(\Omega)$, the $H^2(\Omega)$-norm is not equivalent to the graph norm induced by the operator $\alpha \Delta^2 + \beta \nabla \cdot \nabla (\vec{n}^\top (\mathcal{H} \cdot) \vec{n})$, which would be required for completeness in the natural norm.}

\bibliographystyle{amsplain}
\bibliography{references}
\end{document}